\documentclass[12pt]{article}
\usepackage[centertags]{amsmath}
\usepackage{amsthm,enumerate}
\usepackage[psamsfonts]{amsfonts,amssymb}

\title{Integration by parts on the law of
the reflecting Brownian motion}
\author{{\Large Lorenzo Zambotti} \\
Dipartimento di Matematica
\\ Politecnico di Milano \\ Piazza Leonardo da Vinci 32 
\\ 20133 Milano, Italy\\ \texttt{zambotti@mate.polimi.it}}
\date{}

\numberwithin{equation}{section}

\newtheorem{theorem}{Theorem}[section]
\newtheorem{proposition}[theorem]{Proposition}
\newtheorem{lemma}[theorem]{Lemma}
\newtheorem{corollary}[theorem]{Corollary}

\newcommand{\B}{\dot{B}}

\begin{document}

\maketitle

\begin{abstract}
We prove an integration by parts formula on the law of the
reflecting Brownian motion $X:=|B|$ in the positive half line, 
where $B$ is a standard Brownian motion. In other terms, we consider
a perturbation of $X$ of the form $X^\epsilon = X+\epsilon h$
with $h$ smooth deterministic function and $\epsilon>0$ and we 
differentiate the law of $X^\epsilon$ at $\epsilon=0$. This
infinitesimal perturbation changes drastically the set of zeros 
of $X$ for any $\epsilon>0$. As a consequence, the
formula we obtain contains an infinite dimensional
generalized functional in the sense of Schwartz, defined in
terms of Hida's renormalization of the squared derivative of
$B$ and in terms of the local time of $X$ at $0$.
We also compute the divergence on the
Wiener space of a class of vector fields not taking values 
in the Cameron-Martin space.
\end{abstract}
%


\section{Introduction}

In this paper we want to prove an infinite dimensional 
integration by parts
formula with respect to the law of the reflecting Brownian
motion (RBM) $X_\theta:=|B_\theta-a|$, $\theta\in[0,1]$,
where $B$ is a standard Brownian motion and $a\in{\mathbb
R}$. 

Integration by parts formulae on infinite dimensional
probability measures are an important tool in a number of
topics in Stochastic Analysis. Typically, given a
stochastic process $X$, one considers the law of
an infinitesimal variation
$X^\varepsilon:=X+\varepsilon h$, where $h$
is a process in a suitable class, and one tries to
differentiate the law of $X^\varepsilon$ w.r.t.
$\varepsilon$ at $\varepsilon=0$. 
In most cases one exploits a
quasi-invariance property, i.e. one chooses $h$ in such a 
way that the law of $X^\varepsilon$ is absolutely
continuous w.r.t. the law of $X$: see the monograph 
\cite{ustza}. If this is possible, then the problem is reduced 
to differentiate the density. 

This project has been implemented e.g. for a large class of 
diffusions in ${\mathbb R}^d$ or in Riemannian manifolds,
see e.g. \cite{dri}, \cite{mal} and \cite{hsu}, and for Poisson
measures, see e.g. \cite{akr}. Recently integration by parts for a
class of processes with values in $(0,\infty)$, the Bessel 
bridges of dimension $d\geq 3$, have been computed:
see \cite{za1} and \cite{za2}.

However, the case of processes with a non-trivial
behavior at a boundary remains an open problem. A typical
example of such processes is the RBM, which takes values in
$[0,\infty)$ and has a local time at the boundary $\{0\}$.

In \S 4 of \cite{bis}  J.-M. Bismut developed
a stochastic calculus of variations for the RBM $X=|B-a|$, with the 
aim of studying transition probabilities of boundary
processes associated with diffusions. However
the results of \cite{bis} concern only variations
$X+\varepsilon h$ of $X$ with the crucial
property $\{t:h_t=0\}=\{t:X_t=0\}$. In this case the 
quasi-invariance property holds. Notice that 
$h$ is necessarily a non-deterministic process.

In this paper we consider perturbations $X^\varepsilon= X+\varepsilon h$
of $X=|B-a|$, with $h$ smooth deterministic
function with compact support in $(0,1)$. In this case, the 
approach based on the quasi-invariance fails, since the law
of $X^\varepsilon$ is not absolutely continuous w.r.t. the law of $X$
if $\varepsilon>0$ and $h$ not identically $0$:
see the argument at the end of this introduction.

As a  consequence of the lack of quasi-invariance,
the integration by parts formula we obtain does not
contain only the law of $X$ times suitable densities,
as it is usual in the Malliavin calculus,
see e.g. \cite{mal}, but also an infinite dimensional
generalized functional, in the sense of Schwartz: see Theorem
\ref{main2} below.

This generalized functional is defined in terms of Hida's
square of the white
noise, i.e. a renormalization of the squared derivative of
$B$, defined e.g. in \cite{hkps}, and in terms of the local 
time of $B$ at $0$: see Theorem \ref{main1} below.

It turns out that this problem is closely related with the
computation of the divergence on the Wiener space of a class
of vector fields not taking values in the Cameron-Martin
space. The divergence of vector fields taking values in the 
Cameron-Martin space is typically an $L^p$-variable: see
the monograph \cite{mal}. The divergence we obtain
is not an $L^p$-variable but a generalized functional related with 
the one discussed above: see Theorem \ref{main1.5} below. 

\medskip\noindent
We show now that the law
of $X^\varepsilon$ is not absolutely continuous w.r.t. the law of 
$X=|B-a|$ if $\varepsilon>0$ and $h$ is not identically $0$.
In the case $\min h<0$, with positive probability $\min X^\varepsilon<0$,
while $X\geq 0$ almost surely, so we can suppose $h\geq 0$.
Let $I$ be a non-empty interval where $h>0$ and 
define the set of continuous paths over $[0,1]$:
\[
\Omega^\varepsilon \, := \, \left\{\omega: \min_{\tau\in I}(\omega_\tau-
\varepsilon h_\tau)=0 \right\}.
\]
We claim that ${\mathbb P}(X^\varepsilon\in\Omega^\varepsilon)>0$ while
${\mathbb P}(X\in\Omega^\varepsilon)=0$.

Indeed, $X^\varepsilon\in\Omega^\varepsilon$ if and only if
there exists $\tau\in I$ such that $B_\tau=a$. Since this
event has positive probability, then 
${\mathbb P}(X^\varepsilon\in\Omega^\varepsilon)>0$.
On the other hand:
\[
{\mathbb P}(X\in\Omega^\varepsilon) \, = \, 
{\mathbb P}(B-a\in\Omega^\varepsilon) \, + \, 
{\mathbb P}(a-B\in\Omega^\varepsilon).
\]
By the Girsanov Theorem, the law of $(B_\tau-a-\varepsilon h_\tau:
\tau\in I)$ is absolutely continuous w.r.t. the law of $(B_\tau:
\tau\in I)$, with Radon-Nikodym density $\rho$. In particular:
\[
{\mathbb P}(B-a\in\Omega^\varepsilon) \, = \, 
{\mathbb E}\left[ \rho \ 1_{(\min_I B=0)} \right],
\]
but the r.v. $\min_I B$ has a continuous density, so that
${\mathbb P}(\min_I B=0)=0={\mathbb P}(B-a\in\Omega^\varepsilon)$. Arguing
analogously for ${\mathbb P}(a-B\in\Omega^\varepsilon)$ we obtain
that ${\mathbb P}(X\in\Omega^\varepsilon)=0$.

\section{Main results}

Let $(B_\theta:\theta\in[0,1])$ be a standard Brownian
motion and $C:=\{k:[0,1]\mapsto{\mathbb R}$ continuous,
$k_0=0\}$. We denote by $\mu$ the law of $B$ on $C$: then
$(C,\mu)$ is the classical Wiener space.
We introduce $L:=L^2(0,1)$ with scalar product:
\[
\langle h,k \rangle \, := \, \int_0^1 k_\theta \,
h_\theta \, d\theta, \qquad 
\| h\|^2 \, := \langle h,h \rangle, \qquad h,k\in L.
\]
We consider the following function space on $L$: the set
${\rm Lip}_e(L)$ of $F:L\mapsto{\mathbb R}$ such that:
\[
\exists \ c>0: \quad
|F(h)-F(k)| \, \leq \, e^{c\|h\|} \, \|h-k\|, \qquad h,k\in
L.
\]
Notice that all functions in ${\rm Lip}_e(L)$
are Lipschitz on balls of $L$, with constant growing at most
exponentially with the radius. 

Let $(\rho_\epsilon)_{\epsilon>0}$ be a family of smooth
symmetric mollifiers on ${\mathbb R}$, i.e. 
\[
\rho_\epsilon \, := \, \frac 1\epsilon \rho\left(\frac
\cdot\epsilon\right), \quad \rho\in C^\infty_c(-1,1), \quad
\rho\geq 0, \quad \int_0^1 \rho \, dx \, = \, 1, \quad
\rho(x) \, = \, \rho(-x).
\]
We denote for $\theta\in[0,1]$, $\ell\in C$:
\[
\ell_{\epsilon,\theta} \, = \, (\rho_\epsilon * \ell)_\theta
\, = \, \int_0^1 \rho_\epsilon(\sigma-\theta) \,
\ell_\sigma \, d\sigma,
\]
\[
\stackrel{\textstyle .}{\ell}_{\epsilon,\theta}
\ = \  \ell'_{\epsilon,\theta} \ = \
\frac d{d\theta}
\, \ell_{\epsilon,\theta} \, = \, (-\rho_\epsilon' *
\ell)_\theta. 
\]
With this definition, we denote throughout the paper:
\[
: \! \B_{\epsilon,\theta}^2 \! : \quad \stackrel{\rm def}{=} \quad
\left(\B_{\epsilon,\theta}\right)^2 \ - \ {\mathbb E}\left[
\left(\B_{\epsilon,\theta}\right)^2 \right], \qquad \theta\in[0,1]. 
\]
Here we regularize $B$, we differentiate the regularization
$B_{\epsilon,\cdot}$, we square the derivative and finally we center
this r.v. by subtracting the mean.

Let $(L^a_\theta: \theta\in[0,1])$ denote
the local time of $B$ at $a\in{\mathbb R}$, defined by
the occupation times formula:
\begin{equation}\label{otf}
\int_0^\theta \psi(s,B_s) \, ds \, = \,
\int_{\mathbb R} \int_0^\theta \psi(s,a) \, dL^a_s \, da, \qquad
\theta\in[0,1],
\end{equation}
for all bounded Borel
$\psi:[0,\infty)\times {\mathbb R}\mapsto{\mathbb R}$,
see Chapter VI of \cite{reyo}. Finally, 
let $C_c(0,1)$ denote the space of continuous $h$ with
compact support in $(0,1)$ and $C^2_c(0,1)$ the set of $h\in
C_c(0,1)$ with continuous second derivative. 

Then we can state the first Theorem:
\begin{theorem}\label{main1}
For all $h\in C_c(0,1)$ and $F\in{\rm Lip}_e(L)$,
there exists the limit:
\begin{eqnarray}\label{lim}
& & 
\lim_{\epsilon\to 0} \, {\mathbb E}\left[ F(B) 
\int_0^1 h_\theta \, : \! \B_{\epsilon,\theta}^2 \! :  
\, dL^a_\theta \right]
\\ \nonumber \\ \nonumber & & =: \, 
{\mathbb E}\left[ F(B) \int_0^1 h_\theta 
\, : \! \B_\theta^2 \! : \, dL^a_\theta \right]. 
\end{eqnarray}
\end{theorem}
In the r.h.s. of (\ref{lim}), $: \! \B_\theta^2 \! :$ is the
renormalization of
the square of the derivative of $B$, i.e. Hida's
square of the white noise: since $B$ is not differentiable,
the expression $\B^2$ is not well defined; nevertheless,
subtracting to $\B^2$ a diverging constant, we obtain
convergence to a generalized functional on the Wiener space.
This is made rigorous by the White Noise Analysis, a
generalization to infinite dimension of Schwartz's Theory of
Distributions, see e.g. \cite{hkps}.
However the convergence of the particular functional 
defined by (\ref{lim}) does not seem to be covered by
the existing theorems in the literature, because of the 
integration w.r.t. the local time process. 

Notice that Theorem \ref{main1} defines the r.h.s. of (\ref{lim})
through the limit in the l.h.s.: this can be unsatisfactory and
it seems reasonable to look for a direct way of computing the
functional on $F\in{\rm Lip}_e(L)$: this is done in the last
result of the paper, Corollary \ref{dir} below.
We remark that it is crucial for
the application to the RBM given in Theorem \ref{main2}
below that the limit in (\ref{lim}) exists for a large class
of Lipschitz-continuous functions on $L$, like ${\rm Lip}_e(L)$. 

\medskip\noindent
Before stating the second Theorem, we need a few more notations.
We introduce the Cameron-Martin space
$H^1:=\{h\in C: h'\in L, \ h(0)=0\}$.
We also consider a second function space on $L$:
the set $C^1_e(L)$
of all $F\in {\rm Lip}_e(L)$ with continuous Fr\'echet
differential $\nabla F:L\mapsto L$. Notice that $\nabla F$ satisfies:
\[
\exists \ c>0: \quad
\|\nabla F(h)\| \, \leq \, e^{c\|h\|}, \qquad h\in L.
\]
For any
$\varphi:{\mathbb R}\mapsto{\mathbb R}$ with continuous 
continuous derivative and any smooth deterministic
$h:(0,1)\mapsto{\mathbb R}$ with compact support,
we can define the following vector field over $C$:
\[
{\cal K}:C\mapsto C, \qquad {\cal K}(\omega) \, := \, h\ 
\varphi'(\omega). 
\]
Notice that ${\cal K}$ does not take values in the Cameron-Martin 
space $H^1$, since in general the regularity of $\varphi'(\omega)$ 
is not better than that of $\omega\in C$. Therefore the
divergence of ${\cal K}$ on the Wiener space can not be
computed with the classical theory of the Malliavin calculus,
see \cite{mal}. One of the results of this paper, given in
Theorem \ref{main1.5}, is the computation of this
non-classical divergence.

During the paper we shall consider $\varphi$ in the
class:
\begin{eqnarray*}
{\rm Conv}({\mathbb R}) & := & \Big\{ \
\varphi_1-\varphi_2, \qquad \varphi_i:{\mathbb
R}\mapsto{\mathbb R} \quad {\rm  convex},
\\ & & \qquad \exists \ c>0 \ : \quad
|\varphi_i'(x)| \, \leq \, e^{c|x|}, \quad
\forall x\in{\mathbb R}, \ i=1,2 \ \Big\}.
\end{eqnarray*}
If $F\in C^1_e(L)$, $h\in C_c(0,1)$ and 
$\varphi\in{\rm Conv}({\mathbb R})$, then we can define the directional
derivative of $F$ at $\omega\in C$ along ${\cal K}(\omega)$: 
\[
\partial_{h\, \varphi'(\omega)} F(\omega) \, := \,
\lim_{\epsilon\to 0} \, \frac {F(\omega+\epsilon \, h \, 
\varphi'(\omega)) - F(\omega)}\epsilon.
\]
\begin{theorem}\label{main1.5}
For all
$\varphi\in{\rm Conv}({\mathbb R})$, $h\in C^2_c(0,1)$ and
$F\in C^1_e(L)$ the following integration by parts
formula holds:
\begin{eqnarray}\label{ncm}
{\mathbb E}\left[ \partial_{h\, 
\varphi'(B)} F(B) \right] & = & - \,
{\mathbb E}\left[ F(B) \, \int_0^1 h''_\theta \,
\varphi(B_\theta) \, d\theta \right]
\\ \nonumber \\ \nonumber & &
+ \, \int_{\mathbb R} \, {\mathbb E}\left[ F(B) \int_0^1 
h_\theta \, : \! \B_\theta^2 \! :
\, dL^a_\theta \right] \varphi''(da).
\end{eqnarray}
\end{theorem}

We notice that an infinitesimal transformation along ${\cal 
K}$ does not preserve the absolute-continuity class of
the Wiener measure.
For instance, in the case $\varphi(r)=r^2$, the
infinitesimal transformation along ${\cal K}$ is
$B\mapsto B+\epsilon h\varphi'(B)=B(1+2\epsilon h)$
and it is well known that the laws of $B$ and
$B(1+2\epsilon h)$ are singular if $\epsilon h\ne 0$.
This explains why the r.h.s. of (\ref{ncm}) contains a term, 
the second one, which is not a measure 
but a generalized functional over $C$. We treat the case 
$\varphi(r)=r^2$ and $h\equiv 1$ separately in section 7.

\medskip

We can now turn to the reflecting Brownian motion
$X:=|B-a|$, for some $a\geq 0$. 
For all smooth $f:C\mapsto{\mathbb R}$ and $h\in
C^2_c(0,1)$, by applying (\ref{ncm}) to $F(\omega):=f(|\omega-a|)$
we obtain the following:
\begin{theorem}\label{main2}
We set $X:=|B-a|$ and we denote by $(\ell^0_\theta:
\theta\in[0,1])$ the local time of $X$ at $0$.
Then for all $h\in C^2_c({\mathbb R})$ and $f\in C^1_e(L)$:
\begin{equation}\label{ibp}
{\mathbb E}\left[ \partial_h f(X) \right] \, = \, - \,
{\mathbb E}\left[ f(X) \, \int_0^1 h''_\theta
\, X_\theta \, d\theta \right] 
+ \, {\mathbb E}\left[ f(X) \int_0^1 h_\theta 
\, : \! \B_\theta^2 \! :
\, d\ell_\theta^0 \right].
\end{equation}
\end{theorem}
\noindent
By Tanaka's formula $\ell^0 \equiv 2\, L^a$, see Chapter VI of
\cite{reyo}. Moreover $f(X)=f(|B-a|)$. Therefore the
second term in the r.h.s. of (\ref{ibp}) is defined by
(\ref{lim}).

\medskip\noindent
We give a heuristic argument motivating the result of
Theorem \ref{main1.5}. If $F\in C^1_e(L)$, 
then the classical integration
by parts formula for the Wiener measure states:
\[
{\mathbb E}\left[ \partial_h F(B) \right]  \, = \, 
{\mathbb E}\left[ F(B) \int_0^1 h'_\theta \, dB_\theta \,
\right],
\]
for all deterministic $h\in H^1$, i.e. such that $h'\in
L^2(0,1)$ and $h(0)=0$.

Consider now a process $({\cal K}_\theta(B):\theta\in[0,1])$ 
such that:
\begin{itemize}
\item[1.]  ${\cal K}_\theta = \int_0^\theta \dot{\cal K}_s \, ds$,
with $\dot{\cal K}(B)$ adapted and uniformly bounded.
\item[2.]
there exists a continuous $({\cal Q}_{\theta,\theta'}(\omega):
\theta,\theta'\in[0,1])$ s.t. for all $k\in H^1$:
\[
\left. \frac d{d\varepsilon} \, {\cal K}_\theta(\omega+
\varepsilon k)  \right|_{\varepsilon=0} \, = \, \int_0^1
{\cal Q}_{\theta,\theta'}(\omega) \, k_{\theta'} \, d\theta', \qquad 
\theta\in[0,1], \ \omega\in L.
\]
\end{itemize}
Then the integration by parts formula becomes:
\[
{\mathbb E}\left[ \partial_{{\cal K}(B)} \, F(B) \right] 
\, = \, {\mathbb E}\left[ \left(\int_0^1 \dot{\cal K}_\theta(B) 
\, dB_\theta \, - \, \int_0^1 {\cal Q}_{\theta,\theta}(B) \, 
d\theta \right) F(B) \right].
\]
We set now ${\cal K}_\theta(\omega):=h_\theta \,
\varphi'(\omega_\theta)$, where $h\in C^2_c(0,1)$ and
$\varphi:{\mathbb R}\mapsto{\mathbb R}$ is twice
continuously differentiable with bounded derivatives.
In this case ${\cal K}_\cdot$
is adapted but not a.s. in $H^1$, since $\varphi'(B_\cdot)$ 
has a non-trivial martingale part. Moreover for all $k\in
H^1$: 
\[
\left. \frac d{d\varepsilon} \, {\cal K}_\theta(\omega+ 
\varepsilon k)  \right|_{\varepsilon=0} \, = \, 
h_\theta \, \varphi''(\omega_\theta)\, k_\theta, \qquad
\theta\in[0,1],
\]
so that ${\cal Q}_{\theta,\theta'} = h_\theta \,
\varphi''(\omega_\theta) \, \delta(\theta-\theta')$, where
$\delta$ is the Dirac function. In particular 
${\cal Q}_{\theta,\theta}=h_\theta \,
\varphi''(\omega_\theta) \, \delta(0)$ is
ill-defined, since $\delta(0)=\infty$. However,
arguing formally, we can write:
\[
\int_0^1 {\cal Q}_{\theta,\theta}(B) \, d\theta \, = \,
\int_0^1  h_\theta \,
\varphi''(B_\theta) \, \delta(0) \, d\theta.
\]
Moreover, pretending that $B_{\cdot}$ is
differentiable and $dB_\theta= \ \B_\theta \,  d\theta$,
we obtain:
\begin{eqnarray*}
\int_0^1 \dot{\cal K}_\theta(B) \, dB_\theta & = &
\int_0^1 \frac d{d\theta}\left[h \, \varphi'(B)\right]
\, \B \, d\theta \\ \\ & = &
\int_0^1 h' \, \varphi'(B) \, \B \, d\theta
\, + \int_0^1 h \, \varphi''(B) \, \B^2 \, d\theta.
\end{eqnarray*}
Since $\varphi'(B_\theta) \, \B_\theta = \frac d{d\theta}
\varphi(B_\theta)$, 
integrating by parts over $[0,1]$ in the first
term of this sum, we obtain:
\begin{eqnarray*}
& &
\int_0^1 \dot{\cal K}_\theta(B) \, dB_\theta - \,
\int_0^1 {\cal Q}_{\theta,\theta}(B) \, d\theta
\\ \\ & &  = \,
- \int_0^1 h'' \,  \varphi(B) \, d\theta
\, + \int_0^1 h \, : \! \B^2 \! : \, \varphi''(B) \, d\theta,
\end{eqnarray*}
where $: \! \B^2 \! : \ =  \B^2 - \delta(0)$.
In order to get (\ref{ncm}) we apply the occupation times
formula (\ref{otf}) formally: 
\[
\int_0^1 h \, : \! \B^2 \! : \, \varphi''(B) \, d\theta 
\, = \, \int_{\mathbb R} \left[
\int_0^1 h \, : \! \B^2 \! : \, dL^a_\theta \right] 
\varphi''(da).
\]

\medskip
The paper is organized as follows. In section 3 we prove
that Theorems \ref{main1} and \ref{main1.5} holds for all
$F$ in a suitable space of test functions. In section 4 we
introduce an infinite dimensional Sobolev space on $C$ and
several related functional analytical tools. We prove  
Theorems \ref{main1}, \ref{main1.5} and \ref{main2} in section 5, 
postponing the proof of
the main estimate, given in Lemma \ref{convu}, to section 6.
Finally, in section 7 we discuss the particular case
of quadratic $\varphi$.

\medskip\noindent
We denote by $C_b({\mathbb R})$ the space of
bounded continuous real functions on ${\mathbb R}$ and by
$C^k_b({\mathbb R})$ the set of $f\in C_b({\mathbb R})$ such
the $i$-th derivative of $f$ belongs to $C_b({\mathbb R})$
for all $i=1,\ldots,k$. 

We will use the letter $\kappa$ to denote positive finite
constants whose exact value may change from line to line.

\section{White noise calculus}

In this section we prove that formulae (\ref{lim}) and
(\ref{ncm}) hold for all $F$ in the following
space of test functions over $C$:
\[
{\rm Exp}(C) \, := \, {\rm Span}\{ 
\exp(\langle\cdot,k\rangle):\ k\in C\},
\]
i.e. we prove the following:
\begin{proposition}\label{pr1}
Let $h\in C_c(0,1)$, and $a\in{\mathbb R}$.
Then for all $F\in {\rm Exp}(C)$ 
the limit in (\ref{lim}) exists. 
\end{proposition}
\begin{proposition}\label{pr1.5}
Let $h\in C_c^2(0,1)$ and $\varphi\in
{\rm Conv}({\mathbb R})$.
Then for all $F\in {\rm Exp}(C)$ 
formula (\ref{ncm}) holds. 
\end{proposition}
Propositions \ref{pr1} and \ref{pr1.5} show that Theorems
\ref{main1} and \ref{main1.5} hold for all $F$ in a suitable space of test
functions. The proof of this result is elementary and based
only on the Cameron-Martin theorem and on It\^o's formula.

\medskip

\noindent
We introduce the operator:
\[
Q:L\mapsto L, \qquad
Qk_\theta \, := \, \int_0^1 \theta\wedge\sigma \, k_\sigma
\, d\sigma, \qquad \theta\in[0,1].
\]
The law of $B$ in $L$ is the Gaussian measure with
mean $0$ and covariance operator $Q$, i.e.
\[
{\mathbb E}\left[ e^{\langle B,k\rangle} \right] \, = \,
e^{\frac 12 \langle Qk,k\rangle}, \qquad k\in L.
\]
By the uniqueness of the Laplace transform, we obtain the
following version of the Cameron-Martin formula: 
for all bounded Borel $\Phi:C\mapsto{\mathbb R}$
\begin{equation}\label{cme}
{\mathbb E}\left[ \Phi(B) \, e^{\langle B,k\rangle} \right]
\, = \, e^{\frac 12\langle Qk,k\rangle} \, {\mathbb E}
[\Phi(B+Qk)], \quad k\in C.
\end{equation}
This simple formula is crucial in
White Noise Analysis, in particular in the definition of 
the so called ${\cal S}$-transform: see e.g.
chapter 2 of \cite{hkps}.

We set for $\epsilon<\min\{\theta,1-\theta\}$:
\begin{equation}\label{defce}
c_{\epsilon,\theta} \, := \,
{\mathbb E}\left[\B_{\epsilon,\theta}^2  \right] \, = \, 
\langle Q \rho_\epsilon'(\cdot-\theta),
\rho_\epsilon'(\cdot-\theta) \rangle.
\end{equation}
We also define:
\begin{equation}\label{lambda}
\lambda(\theta,x,y) \, := \, 
x^2 + \, \frac{x \, y}\theta 
\, + \, \frac{y^2-\theta}{4\theta^2}, \qquad
\theta\in(0,1), \ x,y\in {\mathbb R}.
\end{equation}
The proof of Proposition \ref{pr1} is based on the
following: 
\begin{lemma}\label{hi}
For all $\psi\in C_b({\mathbb R})$, $k\in C$,
$K:=Qk$, $\theta\in[\epsilon,1-\epsilon]\subset(0,1)$: 
\begin{equation}\label{otf1}
{\mathbb E}\left[ \psi(B_\theta) \
: \! \B_{\epsilon,\theta}^2 \! : \
e^{\langle B,k\rangle} \right] 
\, = \, e^{\frac12\langle Qk,k\rangle} \, 
{\mathbb E}\left[ \psi(B_\theta+K_\theta) \
\lambda(\theta, K'_{\epsilon,\theta}, B_\theta ) 
\right]. 
\end{equation}
\end{lemma}
\noindent {\bf Proof}. We fix $\theta\in(0,1)$ and set
\[
\ell_\sigma \, := \, 1_{[0,\theta]}(\sigma) \, \frac
\sigma\theta \, + \, 1_{(\theta,1]}(\sigma), \qquad
\beta_\sigma \, := \, B_\sigma - B_\theta \, l_\sigma, \quad
\sigma\in[0,1]. 
\]
Then $\beta$ and $B_\theta$ are independent, i.e.
for all $\Phi:C\mapsto{\mathbb R}$ bounded Borel:
\[
{\mathbb E}[\, \psi(B_\theta)\, \Phi(B)] \, = \,
\int_{\mathbb R} {\cal N}(0,\theta)(dy) \ \psi(y) \,
{\mathbb E}\left[\Phi\left(\beta+y\, \ell\right)\right]. 
\]
Then by (\ref{cme}):
\begin{eqnarray*}
\! & &
{\mathbb E}\left[\, \psi(B_\theta) \,
: \! \B_{\epsilon,\theta}^2 \! : \,
e^{\langle B,k\rangle} \right] \, = \,
\\ \\ & & = \,
e^{\frac12\langle Qk,k\rangle}
{\mathbb E}\left[\, \psi(B_\theta+K_\theta) \,
\left[\left((B+K)'_{\epsilon,\theta}\right)^2
-c_{\epsilon,\theta} \right] \right]
\\ \\ & & = \, e^{\frac12\langle Qk,k\rangle} \!
\int_{\mathbb R} {\cal N}(0,\theta)(dy) \ \psi(y+K_\theta)
\left[ {\mathbb E}\left[
\left((\beta+y\, \ell+K)'_{\epsilon,\theta}\right)^2
\right] -c_{\epsilon,\theta}\right].
\end{eqnarray*}
Since $\theta\in[\epsilon,1-\epsilon]$, we have:  
\[
\ell'_{\epsilon,\theta} \, = \, (\rho_\epsilon *
\ell')_\theta \, = \, \int \rho_\epsilon(\sigma-\theta)
\, \frac 1\theta \,
1_{[0,\theta]}(\sigma) \, d\sigma \, = \, \frac 1{2\theta}.
\]
Then easy computations yield:
\begin{eqnarray*}
& &
{\mathbb E}\left[ \left(
(\beta+y\, \ell+K)'_{\epsilon,\theta} \right)^2 \right]
-c_{\epsilon,\theta}
\\ \\ & & \, = \, \left(
K'_{\epsilon,\theta}\right)^2
\, + \, 2y \, K'_{\epsilon,\theta}
\ell'_{\epsilon,\theta} \, + \, y^2\left(
\ell'_{\epsilon,\theta}\right)^2 \, + \, {\mathbb E}\left[
\left(\beta'_{\epsilon,\theta}\right)^2\right]
-c_{\epsilon,\theta}
\\ \\ & & = \,
\left(K'_{\epsilon,\theta}\right)^2 \,
+ \, \frac y\theta \, K'_{\epsilon,\theta} \, + \,
\frac 1{4\theta^2} \, (y^2-\theta).
\end{eqnarray*}
This yields the thesis. \quad $\square$

\medskip

\noindent
{\bf Proof of Proposition \ref{pr1}.}
Let $h\in C_c(0,1)$. Multiplying (\ref{otf1}) by
$h_\theta$ and integrating in $\theta$ we have:
\begin{eqnarray*}
& &
{\mathbb E}\left[ e^{\langle B,k\rangle}
\int_0^1 h_\theta \, 
: \! \B_{\epsilon,\theta}^2 \! : \,
\psi(B_\theta) \, d\theta \right] 
\\ \\ & & = \, {\mathbb E}\left[ e^{\langle B,k\rangle}
\int_0^1 h_\theta \, \lambda(\theta,K'_{\epsilon,\theta}, 
B_\theta-K_\theta) \, \psi(B_\theta) \, d\theta
\right].
\end{eqnarray*}
By the occupation times formula (\ref{otf}),
this implies for all $a\in{\mathbb R}$: 
\[
{\mathbb E}\left[ e^{\langle B,k\rangle} 
\int_0^1 h_\theta : \! \B_{\epsilon,\theta}^2 \! : \,
dL^a_\theta \right] 
\, = \, {\mathbb E}\left[ 
e^{\langle B,k\rangle} \int_0^1 h_\theta \,
\lambda(\theta,K_{\epsilon,\theta}',a-K_\theta) \,
dL^a_\theta \right].
\]
Since for all $k\in C$ we have $K'_{\epsilon,\theta}\to 
K'_\theta$ as $\epsilon\to 0$, we obtain:
\begin{eqnarray}\label{u7}
& &  \nonumber
{\mathbb E}\left[ e^{\langle B,k\rangle} 
\int_0^1 h_\theta : \! \B_\theta^2 \! : \,
dL^a_\theta \right] \, := \, 
\lim_{\epsilon\to 0}
{\mathbb E}\left[ e^{\langle B,k\rangle} 
\int_0^1 h_\theta : \! \B_{\epsilon,\theta}^2 \! : \,
dL^a_\theta \right] 
\\ \nonumber \\
& & = \, {\mathbb E}\left[ 
e^{\langle B,k\rangle} \int_0^1 h_\theta \,
\lambda(\theta,K_\theta',a-K_\theta) \,
dL^a_\theta \right]. \quad \square
\end{eqnarray}

\medskip\noindent
In Lemma \ref{hi} we have in fact computed the Laplace
transform of the distribution on the Wiener space defined by 
(\ref{lim}): 
\begin{corollary}\label{lap}
For all $a\in{\mathbb R}$, $h\in C_c(0,1)$ and $k\in C$:
\begin{eqnarray*}
& &
{\mathbb E}\left[ e^{\langle B,k\rangle} \int_0^1 h_\theta 
\, : \! \B_\theta^2 \! :
\, dL^a_\theta \right]
\\ \\ & & = \,
e^{\frac12\langle Qk,k\rangle} \, \int_0^1 h_\theta \,
\frac{e^{-(a-K_\theta)^2/2\theta}}{\sqrt{2\pi\theta}} \, 
\lambda(\theta,K_\theta',a-K_\theta) \,
d\theta,
\end{eqnarray*}
where $\lambda$ is defined in (\ref{lambda}).
\end{corollary}

\medskip

\noindent
We turn now to the proof of Proposition \ref{pr1.5}.
For $\Psi_k:=\exp(\langle\cdot,k\rangle)$,
$k\in C$, we have:
\[
\partial_{h \, \varphi'(\omega)}\Psi_k(\omega) \, = \,
\lim_{\epsilon\to 0}
\frac{\Psi_k(\omega+\epsilon \, h\, \varphi'(\omega))
- \Psi_k(\omega)}\epsilon \, = \, \Psi_k(\omega)
\int_0^1 k_\theta h_\theta \, 
\varphi'(\omega_\theta) \, d\theta.
\]
Therefore, by (\ref{cme}),
the l.h.s. of (\ref{ncm}) with $F=\Psi_k$ is equal to: 
\begin{eqnarray}\label{eqt}
{\mathbb E}\left[ \partial_{h \, \varphi'(B)} \Psi_k(B) \right] & = &
{\mathbb E}\left[ \Psi_k(B) \int_0^1 k_\theta \, h_\theta \, 
\varphi'(B_\theta) \, d\theta \right] 
\\ \nonumber \\ \nonumber & = &
e^{\frac 12\langle Qk,k\rangle} \int_0^1 h_\theta \, k_\theta \, 
{\mathbb E}\left[\varphi'(B_\theta+K_\theta)\right] \, d\theta.
\end{eqnarray}
The proof of Proposition \ref{pr1.5} is based on the
following easy application of It\^o's formula.
\begin{lemma}\label{u}
For all $\varphi\in C^2_b({\mathbb R})$, $k\in C$, $K:=Qk$ and
$\theta\in(0,1)$ we have:
\begin{eqnarray}\label{u2}
& & 
k_\theta \, {\mathbb E}\left[
\varphi'(B_\theta+K_\theta) \right] 
\\ \nonumber \\ \nonumber & &
= \, - \, \frac{d^2}{d\theta^2} \, {\mathbb E}\left[
\varphi(B_\theta+K_\theta) \right]
+ \, {\mathbb E}\left[ \varphi''(B_\theta+K_\theta)
\ \lambda(\theta, K'_{\theta}, B_\theta) \right].
\end{eqnarray}
\end{lemma}
\noindent {\bf Proof}. By approximation,
it is enough to consider the case $\varphi\in C^4_b({\mathbb R})$.
By It\^o's formula:
\begin{eqnarray*}
\varphi(B_\theta+K_\theta) & = & \varphi(0) \, + 
\int_0^\theta \varphi'(B_\sigma+K_\sigma) \, (dB_\sigma \, +
\, K_\sigma' \, d\sigma) 
\\ \\ & & + \, \frac 12\int_0^\theta
\varphi''(B_\sigma+K_\sigma) \, d\sigma.
\end{eqnarray*}
Taking expectation and differentiating in $\theta$
we obtain:
\[
\frac d{d\theta} \, {\mathbb E}\left[
\varphi(B_\theta+K_\theta) \right] \, = \, K'_{\theta} \, 
{\mathbb E}\left[ \varphi'(B_\theta+K_\theta)\right]
\, + \, \frac 12 \, {\mathbb E}\left[
\varphi''(B_\theta+K_\theta) \right].
\]
By iteration of this formula we obtain:
\begin{eqnarray*}
& &
\frac{d^2}{d\theta^2} \, {\mathbb E}\left[
\varphi(B_\theta+K_\theta) \right] \, = \,
- \, k_\theta \, {\mathbb E}\left[
\varphi'(B_\theta+K_\theta) \right]
\\ \\ & & + \,
\left( K'_{\theta} \right)^2  {\mathbb E}\left[
\varphi''(B_\theta+K_\theta) \right] 
+ \, K'_{\theta} \,
{\mathbb E}\left[ \varphi'''(B_\theta+K_\theta)
\right] + \, \frac 14 \, {\mathbb E}\left[
\varphi''''(B_\theta+K_\theta)\right].
\end{eqnarray*}
Applying the integration by parts formulae:
\[
\theta \int_{\mathbb R}
\psi'(y+K_\theta) \, {\cal N}(0,\theta)(dy) \, = \, 
\int_{\mathbb R} y\, \psi(y+K_\theta) \, {\cal
N}(0,\theta)(dy),
\]
\[
\theta^2 \int_{\mathbb R}
\psi''(y+K_\theta) \, {\cal N}(0,\theta)(dy) \, = \,
\int_{\mathbb R} (y^2-\theta)\, \psi(y+K_\theta) \, {\cal
N}(0,\theta)(dy),
\]
to $\psi=\varphi''$, we obtain (\ref{u2}). \quad $\square$

\medskip

\noindent
{\bf Proof of Proposition \ref{pr1.5}}. Let $h\in
C^2_c(0,1)$. By a density argument we can reduce to the case
$\varphi\in C^2_b({\mathbb R})$. Multiplying (\ref{u2}) by
$h_\theta$ and integrating in $\theta$
we have, recalling (\ref{lambda}):
\begin{eqnarray}\label{u3}
& & \nonumber
\int_0^1 h_\theta \, k_\theta \, {\mathbb E}\left[
\varphi'(B_\theta+K_\theta) \right] \, d\theta \, = \,
- \, \int_0^1 h''_\theta \ {\mathbb E}\left[
\varphi(B_\theta+K_\theta) \right] \, d\theta 
\\ \nonumber \\  & &
+ \, {\mathbb E}\left[\int_0^1 h_\theta \,
\lambda(\theta,K_\theta',B_\theta) 
\, \varphi''(B_\theta+K_\theta) \, d\theta \right], 
\end{eqnarray}
where $\lambda$ is defined by (\ref{lambda}).
By (\ref{cme}), (\ref{eqt}) and the occupation times
formula (\ref{otf}) this yields:
\begin{eqnarray}\label{u4}
& & \nonumber
{\mathbb E}\left[ \partial_{h \, \varphi'(B)} \Psi_k(B) \right] 
\, = \, - \, \int_0^1 h''_\theta \ {\mathbb E}\left[
\varphi(B_\theta) \, \Psi_k(B) \right] \, d\theta 
\\ \nonumber \\  & & + \, \int_{\mathbb R}
{\mathbb E}\left[\Psi_k(B) \int_0^1 h_\theta \, 
\lambda(\theta,K_\theta',a-K_\theta)
\, dL^a_\theta \right] \varphi''(a) \, da.
\end{eqnarray}
Therefore we conclude by (\ref{u7}). \quad $\square$

\section{Dirichlet forms on the Wiener space}

In this section we introduce infinite dimensional Sobolev spaces
which we need as spaces of test functions. 
Since we consider vector fields ${\cal K}$ taking
values in $C$ or $L$ rather than in the Cameron-Martin
space $H^1$, then the Malliavin derivative is not the correct
notion of gradient and we must introduce a different
differential calculus on $L$.

For $F\in{\rm Exp}(C)$, the usual derivative operator
in the Malliavin calculus is $DF:C\mapsto L$,
defined as follows:
\[
\langle DF(\omega), \ell' \rangle \, := \, \left.
\frac d{d\epsilon} \, F(\omega +\epsilon \ell) \,
\right|_{\epsilon=0}, \qquad \ell\in H^1,
\]
see e.g. \S 1.2 of \cite{nua}. Moreover we have
closability in $L^2(\mu)$ of: 
\[
{\cal D}(F,F) \, := \, \frac 12 \, {\mathbb E}
\left[ \|DF(B)\|^2 \right], \qquad
F\in {\rm Dom}(D) \, = \, {\rm Dom}({\cal D}),
\]
and ${\cal D}$ is a Dirichlet form on the Wiener space. Then
all functions in ${\rm Dom}({\cal D})$ are differentiable in
a weak sense along $H^1$-valued vector fields.

On the other hand we want to study $\partial_{h\,
\varphi'(\omega)} F(\omega)$, see the l.h.s. of (\ref{ncm}),
and in general the regularity of
$\theta\mapsto h_\theta\, \varphi'(\omega_\theta)$ is not
better than that of $\omega\in C$. In particular
the vector field ${\cal K}(\omega) := h\, \varphi'(\omega)$
is not $H^1$-valued and a general $F\in {\rm Dom}({\cal D})$
can not be differentiated along ${\cal K}$. 

For this reason we must consider here a different gradient 
$\nabla F:C\mapsto L=L^2(0,1)$ of $F\in{\rm Exp}(C)$, 
defined by: 
\[
\langle \nabla F(\omega), \ell \rangle \, := \,
\left. \frac d{d\epsilon} \, F(\omega +\epsilon \ell) \,
\right|_{\epsilon=0}, \qquad \ell\in L,
\]
i.e. $\nabla F$ is the Fr\'echet differential of $F$ in $L$. 
Also in this case we have closability in $L^2(\mu)$ of
\[
{\cal E}(F,F) \, := \, \frac 12 \, {\mathbb E}
\left[ \|\nabla F(B)\|^2 \right], \qquad
F\in {\rm Dom}(\nabla) \, = \, {\rm Dom}({\cal E}),
\]
and ${\cal E}$ is a Dirichlet form on the Wiener space.
Comparing the definitions of $DF$ and $\nabla F$ we obtain
$D={\cal P}\nabla$ for all $F\in{\rm Exp}(C)$, where:
\[
{\cal P}:L\mapsto L, \quad {\cal P}\ell_\theta \, := \,
\int_\theta^1 \ell_\tau \, d\tau, \quad \theta\in[0,1].
\]
In particular for some constant $\kappa>0$:
\[
{\cal E}(F,F) \, \geq \, \kappa \, {\cal D}(F,F),
\qquad \forall \, F\in {\rm Dom}({\cal E})
\subset {\rm Dom}({\cal D}).
\]
For a discussion of these infinite dimensional Sobolev
spaces, we refer to \S 9.2.1 for ${\rm Dom}({\cal E})$ 
and to \S 9.3 for ${\rm Dom}({\cal D})$ in \cite{dpz3}. 
We recall in particular that ${\rm Dom}({\cal E})$ also
admits a description in term of the It\^o-Wiener decomposition: see
e.g. Theorem 9.2.12 in \cite{dpz3}.

Now all functions in ${\rm Dom}({\cal E})$ can be
differentiated, at least in a weak sense, along vector
fields taking values in $L$ or $C$, in particular along
${\cal K}(\omega) = h\, \varphi'(\omega)$.
Moreover for $h\in C$ and
$\varphi\in {\rm Conv}({\mathbb R})$, setting:
\[
\Phi_{h,\varphi}=
\Phi:C\mapsto{\mathbb R}, \qquad \Phi(\omega) \, := \,
\langle h,\varphi(\omega) \rangle \, = \, \int_0^1 h_\theta
\, \varphi(\omega_\theta) \, d\theta,
\]
then $\Phi\in{\rm Dom}({\cal E})$ and
$\nabla\Phi(\omega) = h\, \varphi'(\omega)$,
i.e. for all $\omega\in C$: 
\[
\langle\nabla\Phi(\omega), \ell\rangle \, =
\, \int_0^1 h_\theta \, \varphi'(\omega_\theta) \,
\ell_\theta \, d\theta.
\]
Then for all $F\in C^1_e(L)$ the l.h.s. of
(\ref{ncm}) is:
\begin{equation}\label{fr}
{\mathbb E}\left[ \partial_{h\, \varphi'(B)} F(B) \right]
\, = \, {\mathbb E}\left[ \langle\nabla F(B),
h\, \varphi'(B) \rangle\right] \, = \,
2 \, {\cal E}(F, \Phi_{h,\varphi}).
\end{equation}

We recall now that the semigroup $(P^{\cal D}_t: t\geq 0)$
in $L^2(\mu)$ associated with ${\cal D}$ is given by the
Mehler formula: 
\[
P^{\cal D}_t F(z) \, = \, \int F(y) \ {\cal N}
\left(e^{-t/2} \, z, (1-e^{-t}) \, Q \right)(dy), \qquad
z\in C, \ F\in L^2(\mu),
\]
where ${\cal N}(a,{\cal Q})$ denotes the Gaussian measure
over $L$ with mean $a\in L$ and covariance operator ${\cal
Q}:L\mapsto L$.
This semigroup is a basic tool in the Malliavin calculus: 
see e.g. Chapters 1-2 in \cite{mal} and \S 1.4-1.5 in \cite{nua}.

Since in this paper we work with $\nabla$ 
rather than with $D$, a crucial role is played by the
transition semigroup $(P_t: t\geq 0)$ in $L^2(\mu)$
associated with ${\cal E}$, given by:
\[
P_t F(z) \, = \, \int F(y) \ {\cal N}
\left(e^{tA} \, z,  \, Q_t \right)(dy), \qquad
z\in C, \ F\in L^2(\mu),
\]
where $(e^{tA}: t\geq 0)$ is the semigroup in $L$ generated
by the operator:
\[
D(A) \, := \, \{h\in C: h''\in L, \ h(0)=h'(1)=0\}, \qquad
Ah \, := \, \frac 12 \, h'',
\]
and we set:
\begin{equation}\label{qmt}
Q_t \, := \, \int_0^t e^{2sA} \, ds \, = \,
\frac{I-e^{2tA}}{-2A}, \qquad t\in[0,\infty].
\end{equation}
Notice in particular that: 
\begin{equation}\label{qinf}
Q_\infty \, = \, (-2A)^{-1} \, = \, Q. 
\end{equation}
The second equality of (\ref{qinf}) says that $Q$ and $-2A$ are inverse
of one another and can be verified by an explicit computation.

The operators $(P^{\cal D}_t: t\geq 0)$ and $(P_t: t\geq 0)$ are 
two different examples of Ornstein-Uhlenbeck semigroups: we refer 
to Chapters 6 and 10 in \cite{dpz3}. For a more detailed
description of $(P_t: t\geq 0)$ see section 6 below. 

Two important properties of ${\rm Dom}({\cal E})$ are stated in
the following:
\begin{lemma}\label{lips}
The space ${\rm Lip}_e(L)$ is contained in ${\rm Dom}({\cal
E})$. The space ${\rm Exp}(C)$ is dense in ${\rm
Dom}({\cal E})$.
\end{lemma}
\noindent{\bf Proof}. We recall that $F\in {\rm Dom}({\cal
E})$ if and only if $\sup_{t>0}{\cal E}(P_tF,P_tF) < \infty$. Now: 
\begin{eqnarray*}& &
|P_t F(z_1) - P_tF(z_2)| \, \leq \, 
\int |F(y+e^{tA}z_1)-F(y+e^{tA}z_2)| \, {\cal N}(0,Q_t)(dy) 
\\ \\ & & \leq \, \int e^{c(\|y\|+\|z_1\|)} \, \|z_1 - z_2\| 
\, {\cal N}(0,Q_t)(dy) \, \leq \, \kappa \,
e^{c\|z_1\|} \, \|z_1 - z_2 \|,
\end{eqnarray*}
so that $\|\nabla P_t F(z)\| \leq e^{c\|z\|}$ for all $z\in C$ and 
we obtain the first claim. For the second one, we refer to
\S 9.2.1 of \cite{dpz3}. \quad $\square$

\section{Proof of the main results}

We want to use the tools
introduced in the previous section to prove Theorems
\ref{main1}, \ref{main1.5} and \ref{main2}.

In Propositions
\ref{pr1} and \ref{pr1.5} we have proved that (\ref{lim}) 
and (\ref{ncm}) hold for all $F\in{\rm
Exp}(C)$. This space is dense in the topology of the
Sobolev space ${\rm Dom}({\cal E})$, introduced in the
previous section. An a priori estimate, given in Lemma
\ref{convu}, and a density argument allow to extend
(\ref{lim}) and (\ref{ncm}) to much larger spaces of
test functions and to prove Theorems \ref{main1} and 
\ref{main1.5} and also Theorem \ref{main2} as a corollary. 
In particular, in this section we prove:
\begin{proposition}\label{pr3}
Let $h\in C_c(0,1)$ and $a\in{\mathbb R}$. Then the limit in 
(\ref{lim}) exists for all $F\in{\rm Lip}_e(L)$. 
\end{proposition}
\begin{proposition}\label{pr4}
For all $h\in
h^2_c(0,1)$, $\varphi\in{\rm Conv}({\mathbb R})$ and
$F\in{\rm Lip}_e(L)$:
\begin{eqnarray}\label{ncm3}
{\mathbb E}\left[ \langle\nabla F(B), h\, 
\varphi'(B) \rangle \right] & = & - \,
{\mathbb E}\left[ F(B) \, \int_0^1 h''_\theta \,
\varphi(B_\theta) \, d\theta \right]
\\ \nonumber \\ \nonumber & &
+ \, \int_{\mathbb R} \, {\mathbb E}\left[ F(B) \int_0^1 
h_\theta \, : \! \B_\theta^2 \! : \, dL^a_\theta \right]
\varphi''(da). 
\end{eqnarray}
\end{proposition}

\noindent
Proposition \ref{pr3} proves Theorem
\ref{main1}. Theorem \ref{main1.5} follows by Proposition
\ref{pr4}  and formula (\ref{fr}),
recalling that $C^1_e(L)\subset{\rm Lip}_e(L)$. At the end
of the section, we derive Theorem \ref{main2} from Proposition
\ref{pr4}. We also recall 
that $\nabla F$ is well defined, since by Lemma \ref{lips}:
$F\in{\rm Lip}_e(L)\subset{\rm Dom}({\cal E})=
{\rm Dom}(\nabla)$.

\medskip
\noindent
We recall that $\mu$ denotes the Wiener measure,
law of $B$, i.e. for all bounded Borel $F:C\mapsto{\mathbb R}$:
\[
\mu(F) \, = \, \int F \, d\mu \, = \, 
{\mathbb E} [F(B)].
\]
By Proposition 10.5.2 of \cite{dpz3},
${\cal E}$ satisfies the Poincar\'e inequality: 
\[
\int \left(F-\mu(F)\right)^2 \, d\mu \, \leq \,
\frac 1{\lambda_1} \  {\cal E}(F,F), \qquad F\in
{\rm Dom}({\cal E}),
\]
where $\lambda_1=\pi^2/4$, see (\ref{eig}) below. 
Since $(P_t: t\geq 0)$
is the semigroup in $L^2(\mu)$ associated with ${\cal E}$, 
the Poincar\'e inequality implies the exponential
convergence of $P_tF$ to $\mu(F)$ in $L^2(\mu)$:
\begin{equation}\label{expc}
\|P_t F-\mu(F)\|^2_{L^2(\mu)} \, \leq \, e^{-2t/\lambda_1} \,
\|F\|^2_{L^2(\mu)}, \qquad t\geq 0, \ F\in L^2(\mu).
\end{equation}
In particular for all $G\in L^2(\mu)$: 
\[
{\cal R}G \, := \, \int_0^\infty \left(P_t G - 
\mu(G) \right) \, dt \, \in \, {\rm Dom}({\cal E}),
\]
and for all $F\in {\rm Dom}({\cal E})$: 
\[
{\mathbb E}[ F(B) \, G(B) ] \, = \,
{\mathbb E}[ F(B)] \, {\mathbb E}[ G(B) ] \, + \,
{\cal E}(F,{\cal R}G).
\]
Let now $h\in C_c(0,1)$ and $a\in{\mathbb R}$.
For all $\epsilon>0$ we define
$G_{\epsilon,a}\in L^2(\mu)$
\begin{equation}\label{f-e}
G_{\epsilon,a}(B) \, := \, \int_0^1 h_\theta \,
: \! \B_{\epsilon,\theta}^2 \! : \, dL^a_\theta, \qquad 
G_\epsilon \, := \, G_{\epsilon,0}.
\end{equation}
Then (\ref{lim}) is equivalent to
the existence of the limit as $\epsilon\to 0$ of:
\begin{equation}\label{calr}
{\mathbb E}\left[F(B) \, G_{\epsilon,a}(B)\right] \, = \,
{\mathbb E}[F(B)] \,{\mathbb E}[ G_{\epsilon,a}(B) ] \, + \,
{\cal E}(F,{\cal R}G_{\epsilon,a})
\end{equation}
for all $F\in {\rm Lip}_e(L)$. The main tool in the proof of
Propositions \ref{pr3} and \ref{pr4} is the following estimate:
\begin{lemma}\label{convu} If $h\in C_c(0,1)$ then there exists
a constant $\kappa>0$ such that:
\begin{equation}\label{dec}
\|P_t G_\epsilon\|^2_{L^2(\mu)} \, \leq \, \kappa \, \frac{
1+|\ln t|^6}{t^{3/4}}, \qquad t\in(0,1], \ \epsilon>0.
\end{equation}
\end{lemma}
\noindent
The proof of Lemma \ref{convu} is postponed to section 6.
As a consequence of Lemma \ref{convu} we have the following: 

\begin{proposition}\label{pr2}
Let $h\in C_c(0,1)$ and $a=0$. Then the limit in (\ref{calr}) 
exists for all $F\in {\rm Dom}({\cal E})$. 
\end{proposition}
\noindent
{\bf Proof.}
By Lemma \ref{hi} for $k=0$ and $\psi\in C_b({\mathbb R})$
we have:
\[
{\mathbb E}\left[ \int_0^1 h_\theta \, 
: \! \B_{\epsilon,\theta}^2 \! : \,
\psi(B_\theta) \, d\theta \right] 
\, = \, {\mathbb E}\left[ \int_0^1 h_\theta \,
\frac{B^2_\theta-\theta}{4\theta^2} \,
\psi(B_\theta) \, d\theta \right].
\]
By the occupation times formula (\ref{otf}) we obtain for all
$\psi\in C_b({\mathbb R})$:
\[
\int_{\mathbb R} {\mathbb E}[G_{\epsilon,a}(B)] \, \psi(a)
\, da \, = \int_0^1 h_\theta \,
\int_{\mathbb R} \frac{a^2-\theta}{4\theta^2} \,
\frac{e^{-a^2/2\theta}}{\sqrt{2\pi\theta}} \,  
\psi(a) \, da \, d\theta.
\]
In particular:
\[
{\mathbb E}[G_{\epsilon,a}(B)] \, = \, \int_0^1 h_\theta \,
\frac{a^2-\theta}{4\theta^2} \,
\frac{e^{-a^2/2\theta}}{\sqrt{2\pi\theta}} \, d\theta,
\]
which does not depend on $\epsilon$. Therefore, by (\ref{calr}) 
the existence of the limit in (\ref{lim}) with $a=0$ for all 
$F\in {\rm Dom}({\cal E})$ is equivalent to the weak convergence of
${\cal R}G_\epsilon$ in ${\rm Dom}({\cal E})$. 
Now, by Proposition \ref{pr1}, the limit in (\ref{lim}) with $a=0$ exists for 
all $F\in {\rm Exp}(C)$, which is dense in ${\rm Dom}({\cal E})$.
Therefore, if we can prove that:
\begin{equation}\label{esti}
\sup_{\epsilon>0} \ {\cal E}({\cal R}G_\epsilon,
{\cal R}G_\epsilon) \, < \, \infty,
\end{equation}
then we conclude. Indeed, for any $F\in {\rm Dom}({\cal E})$ we can find
a sequence $(F_n)_n\subset {\rm Exp}(C)$ converging to $F$ in ${\rm Dom}({\cal E})$.
Write:
\begin{eqnarray*}
& &
|{\cal E}(F,G_\epsilon-G_\delta)| \, \leq \,
|{\cal E}(F_n,G_\epsilon-G_\delta)| \, + \, |{\cal E}(F-F_n,G_\epsilon-G_\delta)|.
\end{eqnarray*}
By (\ref{esti}) we can make the second term arbitrarily small for some $n$ 
big enough but fixed, uniformly in $\epsilon,\delta>0$. Then by Proposition
\ref{pr1} we can make the first term arbitrarily small as $\epsilon,\delta\to 0$.

For the proof of (\ref{esti}), we recall the following formula: 
\begin{eqnarray*}
& &
{\cal E}({\cal R}G_\epsilon,{\cal R}G_\epsilon) \, = \,
\int {\cal R}G_\epsilon \,
(G_\epsilon-\mu(G_\epsilon)) \, d\mu
\\ \\ & & = \,
\int_0^\infty \int (P_tG_\epsilon-\mu(G_\epsilon)) 
\, (G_\epsilon-\mu(G_\epsilon)) \, d\mu \, dt
\\ \\ & & = \, \int_0^\infty 
\|P_{t/2} G_\epsilon-\mu(G_\epsilon)\|^2_{L^2(\mu)} \, dt.
\end{eqnarray*}
Moreover by (\ref{expc}) and (\ref{dec}),
since $P_{1+t} = P_t P_1$, $t\geq 0$: 
\[
\|P_{1+t} G_\epsilon-\mu(G_\epsilon)\|^2_{L^2(\mu)}
\, \leq \, e^{-2t/\lambda_1} \, \|P_1G_\epsilon\|^2_{L^2(\mu)} 
\, \leq \, \kappa \, e^{-2t/\lambda_1}.
\]
Therefore (\ref{esti}) follows from:
\begin{eqnarray*}
{\cal E}({\cal R}G_\epsilon,{\cal R}G_\epsilon) & \leq &
\int_0^1 \|P_{t/2} G_\epsilon\|^2_{L^2(\mu)} \, dt \, + \,
\int_1^\infty \|P_{t/2}
G_\epsilon-\mu(G_\epsilon)\|^2_{L^2(\mu)} \, dt 
\\ \\ & \leq & \kappa \int_0^1 \frac{1+|\ln t|^6}{t^{3/4}} \, dt
\, + \, \kappa \int_1^\infty e^{-2t/\lambda_1} \, dt \, < \, \infty. \quad
\square
\end{eqnarray*}

\noindent
We can now apply the results of Propositions 
\ref{pr1}, \ref{pr1.5} and \ref{pr2} to prove Propositions 
\ref{pr3} and \ref{pr4}.

\noindent{\bf Proof of Proposition \ref{pr3}.}
We fix $\delta\in(0,1/2)$ such that ${\rm
supp}(h)\subset[\delta,1-\delta]$ and we consider
$\epsilon\in(0,\delta/2)$.
By Proposition \ref{pr2}, (\ref{lim}) holds for $a=0$ and
for all $F\in {\rm Dom}({\cal E})$.

Let $\ell:[0,1]\mapsto{\mathbb R}$ be of class $C^2$
such that $\ell_0=0$ and $\ell_\theta=1$ for all
$\theta\in[\delta/2,1]$. By the Cameron-Martin theorem 
we have the following formula:
\begin{equation}\label{ell}
{\mathbb E}\left[ F(B) \right] \, = \, {\mathbb E}\left[ 
F(B+a\ell) \, \exp\left( a\langle \ell'', B \rangle
- c(\ell,a)\right) \right],
\end{equation}
where $c(\ell,a) \, := \, a^2\|\ell'\|^2/2$.
If $G_{\epsilon,a}$ is defined as in (\ref{f-e}), then almost surely:
\[
G_{\epsilon,a}(B+a\ell) \, = \, G_{\epsilon,a}(B+a) \, = \, 
G_{\epsilon,0}(B) \, = \, G_\epsilon(B),
\]
where the first equality holds because $h$ vanishes where
$\ell\ne 1$ and the second one because the local time of $B+a$
at $a$ is equal to the local time of $B$ at $0$.
Let now $F$ be in ${\rm Lip}_e(L)$. Then by (\ref{ell}):
\begin{equation}\label{ell2}
{\mathbb E}\left[ F(B) \, G_{\epsilon,a}(B) \right] \, = \,
{\mathbb E}\left[ F_a(B)\, 
G_\epsilon(B) \right], 
\end{equation}
\[
{\rm where}: \qquad
F_a(z) \, := \, F(z+a\ell)\, \exp(a\langle \ell'', z \rangle
- c(\ell,a)), \quad z\in C.
\]
Now, $F_a\in {\rm Dom}({\cal E})$, so that, by Proposition \ref{pr2},
${\mathbb E}[ F(B) \, G_{\epsilon,a}(B) ]$ converges as
$\epsilon\to 0$ and (\ref{lim}) is proven for all
$a\in{\mathbb R}$. \quad $\square$

\medskip

\noindent{\bf Proof of Proposition \ref{pr4}.}
We consider first the case:
\[
\varphi(x) \, := \, |x-a| \ \Longrightarrow \
\varphi'(x) \, = \, {\rm sign}(x-a), \quad
\varphi''(dx) \, = \, 2\, \delta_a(dx),
\]
for some $a\in{\mathbb R}$, where $\delta_a$ is the Dirac
mass at $a$ and
\[
{\rm sign}:{\mathbb R}\mapsto\{0,1\}, \qquad
{\rm sign}(x) \, := \,
1_{(0,\infty)}(x) \, - \, 1_{(-\infty,0]}(x).
\]
In this case, (\ref{ncm3}) becomes:
\begin{eqnarray}\label{ncmd}
{\mathbb E}\left[ \langle\nabla F(B), h\, 
{\rm sign}(B-a) \rangle \right] & = & - \,
{\mathbb E}\left[ F(B) \, \int_0^1 h''_\theta \,
|B_\theta-a| \, d\theta \right]
\\ \nonumber \\ \nonumber & &
+ \, 2 \, {\mathbb E}\left[ F(B) \int_0^1 h_\theta 
\, : \! \B_\theta^2 \! : \, dL^a_\theta \right].
\end{eqnarray}

Consider first the case $a=0$.
By Proposition \ref{pr2}, the r.h.s. of (\ref{ncmd})
defines a bounded linear functional on ${\rm Dom}({\cal E})$.
Moreover, by Proposition \ref{pr1}, (\ref{ncmd}) holds 
for all $F\in{\rm Exp}(C)$. Since
both sides of (\ref{ncmd}) are bounded 
linear functionals on ${\rm Dom}({\cal E})$, coinciding on the dense 
subset ${\rm Exp}(C)$, they coincide on ${\rm Dom}({\cal E})$.
Therefore (\ref{ncmd}) is proven for $a=0$.

Let $\ell$ and $F_a$ be the functions introduced in the proof of
Proposition \ref{pr3}. By 
(\ref{ell2}) and by (\ref{ncmd}) with $a=0$ we obtain: 
\begin{eqnarray*}
& &
\lim_{\epsilon\to 0} \, 2 \, {\mathbb E}\left[ F(B) \, 
G_{\epsilon,a}(B) \right] \, = \, 
\lim_{\epsilon\to 0} \, 2 \, 
{\mathbb E}\left[ F_a(B)\, G_\epsilon(B) \right] 
\\ \\ & & = \, {\mathbb E}\left[ \langle\nabla F_a(B), h\, 
{\rm sign}(B) \rangle \right] + \,
{\mathbb E}\left[ F_a(B) \, \int_0^1 h''_\theta \,
|B_\theta| \, d\theta \right] 
\\ \\ & & = \, {\mathbb E}\left[ \langle\nabla F(B),
h\, {\rm sign}(B-a) \rangle \right] + \,
{\mathbb E}\left[ F(B) \, \int_0^1 h''_\theta \,
|B_\theta-a| \, d\theta \right].
\end{eqnarray*}
Therefore (\ref{ncmd}) is proven for all $a\in{\mathbb R}$
and $F\in{\rm Lip}_e(L)$.
Let now $\varphi\in C^2_c({\mathbb R})$. Multiplying
(\ref{ncmd}) by $\varphi''(a)$ and integrating in $da$ we
obtain (\ref{ncm3}) and: 
\begin{eqnarray*}
& &
\left|\int_{\mathbb R} \, {\mathbb E}\left[ F(B) \int_0^1 
h_\theta \, : \! \B_\theta^2 \! : \, dL^a_\theta \right] 
\varphi''(a) \, da\right|
\\ \\ & &
\leq \, \kappa \int_0^1 {\mathbb E}\left[ e^{c\|B\|} \left( 
|\varphi(B_\theta)| + |\varphi'(B_\theta)| \right) \right] 
d\theta.
\end{eqnarray*}
Therefore, by a density argument (\ref{ncm3}) holds for all 
$\varphi\in {\rm Conv}({\mathbb R})$. \quad $\square$

\medskip\noindent
{\bf Proof of Theorem \ref{main2}}.
We start by recalling that, by Tanaka's formula,
$\ell^0 \equiv 2 \, L^a$, where
$\ell^0$ is the local time process of $X=|B-a|$ at $0$ and
$L^a$ is the local time process of $B$ at $a$.

Fix $h\in C^2_c(0,1)$ and
$f\in C^1_e(L)$. Setting $F(z):=f(|z-a|)$, $z\in L$, then
clearly $F\in{\rm Lip}_e(L)$. By Lemma \ref{lips}, $F\in{\rm
Dom}({\cal E})$ and by the chain rule:
\[
\langle\nabla F(z),h\rangle \, = \, \langle\nabla f(|z-a|),
h\, {\rm sign}(z-a) \rangle, \qquad \mu-{\rm a.e.} \ z.
\]
In particular for $\mu$-a.e. $z$:
\[
\langle\nabla F(z),h\, {\rm sign}(z-a) \rangle \, = \,
\langle\nabla f(|z-a|), h \rangle,
\]
since $[{\rm sign}(z-a)]^2\equiv 1$. Therefore,
formula (\ref{ncmd}) applied to $F(z):=f(|z-a|)$ and
$\varphi(x)=|x-a|$, $z\in L$, $x\in{\mathbb R}$, yields
(\ref{ibp}). \quad $\square$

\section{The main estimate}

In this section we prove Lemma \ref{convu}. We recall that
$G_\epsilon$ is the sum of two diverging terms. Applying
$P_t$ to $G_\epsilon$ we have a regularization effect:
indeed, we write $P_t G_\epsilon$ as a sum of terms, which
after some cancelations converge as $\epsilon$ tends to $0$.
This compensation of infinities requires a careful study
of each term.

\medskip \noindent
We start with a more detailed description of the semigroup
$(P_t:t\geq 0)$ of the Dirichlet Form ${\cal E}$ in
$L^2(\mu)$, defined in section 4.
We introduce first the Green function
$(g_t(\theta,\theta'):t>0, \theta,\theta'\in[0,1])$ of the
heat equation associated with $A$, i.e. solution of 
\[
\frac{\partial g}{\partial t} \, = \,
\frac 12 \frac{\partial^2 g}{\partial \theta^2}
\]
with boundary and initial conditions:
\[
g_t(0,\theta') \, = \, \frac
{\partial g_t}{\partial\theta}(1,\theta')=0,
\qquad
g_0(\theta,\theta')\, = \, \delta_\theta(d\theta'),
\]
where $\delta_\theta$ is the Dirac
mass at $\theta$. Then we set for all $z\in C$:
\begin{equation}\label{defzv}
z(t,\theta) \, := \int_0^1 g_t(\theta,\theta') \,
z_{\theta'} \, d\theta', 
\quad v(t,\theta) \, := \int_0^t \int_0^1
g_{t-s}(\theta,\theta') \, W(d\theta',ds),
\end{equation}
\begin{equation}\label{defu}
u(t,\theta) \, := \, z(t,\theta) \, + \, v(t,\theta),
\qquad U_t(z) \, := \, u(t,\cdot)\in C,
\end{equation}
where $(W(\theta',s): \theta'\in[0,1], s\geq 0)$ is a Brownian 
sheet. Then:
\[
P_tF(z) \, = \, {\mathbb E}\left[ F(U_t(z)) \right], \qquad
t\geq 0, \ z\in C, \ F\in L^2(\mu).
\]
Although this is not needed in this paper, we remark 
that $(u(t,\theta): t\geq 0, \theta\in[0,1])$ is the unique
solution of the Stochastic Partial
Differential Equation driven by space-time white noise:
\[
\left\{ \begin{array}{ll}
{\displaystyle
\frac{\partial u}{\partial t}=\frac 12
\frac{\partial^2 u}{\partial \theta^2}
 + \frac{\partial^2 W}{\partial t\partial\theta} 
}
\\ \\
{\displaystyle u(t,0)=\frac
{\partial u}{\partial\theta}(t,1)=0 }
\\ \\
u(0,\theta)=z_\theta,
\end{array} \right.
\]
see \cite{wa}.

Notice that $(z(t,\theta): t\geq 0,\theta\in[0,1])$ is a 
deterministic continuous function and $(v(t,\theta):
t\geq 0, \theta\in[0,1])$ is a centered continuous
Gaussian process. A crucial role is played by the function: 
\begin{equation}\label{covv}
q_t(\theta,\theta') \, := \,
{\mathbb E}\left[v(t,\theta) \, v(t,\theta') \right] \, =
\, \int_0^t g_{2s}(\theta,\theta') \, ds, \qquad
q_t(\theta) \, := \, q_t(\theta,\theta),
\end{equation}
for $\theta,\theta'\in[0,1]$, $t\geq 0$.
Notice that for all $\ell\in L$:
\[
\int_0^1 q_t(\theta,\theta') \, \ell_{\theta'} \, d\theta'
\, = \, \int_0^t e^{2sA}\ell_\theta \, ds \, = \,
Q_t\ell_\theta,
\]
where $Q_t$ is defined in (\ref{qmt}). By (\ref{qinf}) above,
$Q_\infty=Q$, i.e.
\begin{equation}\label{qi}
q_\infty(\theta,\theta') \, := \, \lim_{t\nearrow \infty}
q_t(\theta,\theta') \, = \, \theta\wedge\theta', \qquad
q_\infty(\theta) \, := \, q_\infty(\theta,\theta) \, = \,
\theta.
\end{equation}
We also set:
\begin{equation}\label{defqi}
q^t(\theta,\theta') \, := \, [q_\infty - q_t]
(\theta,\theta') \, = \,
\int_t^\infty g_{2s}(\theta,\theta') \, ds,
\qquad q^t(\theta) \, := \, q^t(\theta,\theta).
\end{equation}

We denote by $\gamma_t(\theta-\theta')$ the density of the 
Gaussian measure ${\cal N}(\theta,t)(d\theta')$ over
${\mathbb R}$ with mean $\theta$ and variance $t$.
Then $g-\gamma$ is smooth over $[0,\infty)\times(0,1)
\times(0,1)$. In particular for all $\delta\in(0,1/2)$ 
there exists a constant $\kappa_\delta>0$ such that
for all $t\in[0,1]$, $\theta\in[\delta,1-\delta]$:
\begin{equation}\label{estq}
q_t(\theta) \, = \, \int_0^t \frac{ds}{\sqrt{4\pi s}}
\, + \, \int_0^t (g_{2s}(\theta,\theta) - \gamma_{2s}(0))
\, ds \, \geq \, \kappa_\delta \, t^{1/2}.
\end{equation}
Finally, we introduce the complete orthonormal system of $L$:
\[
e_i(\theta) \, := \, 2^{1/2} \, \sin\left(\sqrt{\lambda_i} \, \theta\right), \quad
\theta\in[0,1], \qquad \lambda_i \, := \, \frac{\pi^2}4 \, (2i-1)^2, 
\]
$i=1,2,\ldots$.
Then $(e_i)_i$ is a system of eigenvectors of $Q$, $A$ and $e^{tA}$:
\begin{equation}\label{eig}
Q \, e_i \, = \, \frac 1{\lambda_i} \, e_i, \qquad
A \, e_i \, = \, - \, \frac{\lambda_i} 2 \, e_i, \qquad
e^{tA} \, e_i \, = \, e^{-t\lambda_i/2} \, e_i.
\end{equation}
In particular:
\begin{equation}\label{cos}
q_t(\theta,\theta') \, = \, \sum_{i=1}^\infty \frac{1-e^{-\lambda_i t}}
{\lambda_i} \, e_i(\theta) \, e_i(\theta'), \qquad t\in[0,\infty], \ \theta,
\theta'\in[0,1].
\end{equation}

\medskip
\noindent{\bf Proof of Lemma \ref{convu}.} 
We fix $\delta\in(0,1/2)$ such that 
${\rm supp}(h)\subseteq[\delta,1-\delta]$ and we consider 
$\epsilon\in(0,\delta)$. Recalling (\ref{defzv}) and
(\ref{defu}), we set
\[
v_\epsilon(t,\cdot) \, := \, \rho_\epsilon* v(t,\cdot),
\quad 
z_\epsilon(t,\cdot) \, := \, \rho_\epsilon* z(t,\cdot),
\quad u_\epsilon \, := \, z_\epsilon+v_\epsilon.
\] 
We denote the partial derivative w.r.t. $\theta$ by
$\partial_\theta$. 

\medskip\noindent
{\bf An explicit formula for $P_tG_\epsilon$.}
By the definition (\ref{f-e}) of $G_{\epsilon,a}$
and by the occupation times formula, for $\mu$-a.e. $\omega$: 
\[
\int_{\mathbb R} G_{\epsilon,a}(\omega)\, \psi(a) \, da \, =
\, \int_0^1 h_\theta \, \left( (\omega'_{\epsilon,\theta})^2
- c_{\epsilon,\theta} \right) \, \psi(\omega_\theta) \,
d\theta,
\]
for any $\psi\in C_b({\mathbb R})$. By Fubini's theorem:
\begin{eqnarray}\label{cc}
\int_{\mathbb R} P_tG_{\epsilon,a}(z) \, \psi(a) \, da & = & 
P_t\left[\int_{\mathbb R} G_{\epsilon,a}\, \psi(a) \, da
\right](z)
\\ \nonumber\\ \nonumber
& = & \int_0^1 h_\theta \,
{\mathbb E} \left[ \psi(u(t,\theta)) 
\left( ( \partial_\theta u_\epsilon(t,\theta))^2
-c_{\epsilon,\theta}\right)  \right] \, d\theta.
\end{eqnarray}
As in the proof of Lemma \ref{hi} we set 
for fixed $t>0$ and $\theta\in(0,1)$:
\[
\ell_\sigma \, := \,
\frac{q_t(\sigma,\theta)}{q_t(\theta)}, \qquad
\hat{v}(t,\sigma) \, := \, v(t,\sigma) - v(t,\theta) \,
\ell_\sigma, \qquad \sigma\in(0,1).
\]
Then the covariance between the two Gaussian variables
$\hat{v}(t,\cdot)$ and $v(t,\theta)$ is zero, so that  
$\hat{v}(t,\cdot)$ and $v(t,\theta)$ are independent.
Denoting $\overline z:=z(t,\theta)$ and
$\overline q:=q_t(\theta)$ we obtain:  
\begin{eqnarray}\label{ccc}
& &
{\mathbb E} \left[ \psi(u(t,\theta)) 
\left[ ( \partial_\theta u_\epsilon(t,\theta))^2
-c_{\epsilon,\theta}\right]  \right]
\\ \nonumber \\ \nonumber & & = 
\int_{\mathbb R} {\cal N}(0,\overline q)(dy) \
\psi(y+\overline z) \, {\mathbb E}\left[\left(\partial_\theta
u_\epsilon(t,\theta) +(y-v(t,\theta)) \, \ell'_{\epsilon,\theta}
\right)^2-c_{\epsilon,\theta}\right]
\\ \nonumber \\ \nonumber & & = 
\int_{\mathbb R} {\cal N}(0,\overline q)(dy) \,
\psi(y+\overline z) \, \left[
\left(\partial_\theta z_\epsilon(t,\theta)
 + \, y \, \ell'_{\epsilon,\theta} \right)^2
\, - \, \overline q
\left(\ell'_{\epsilon,\theta}\right)^2
\, - \, c^t_{\epsilon,\theta} \right]
\end{eqnarray}
where, recalling (\ref{qmt}) and setting 
$Q^t := e^{tA}Qe^{tA} = Q-Q_t$, by (\ref{defce}): 
\begin{equation}\label{defcet}
c^t_{\epsilon,\theta} \, := \, c_{\epsilon,\theta} - \, 
{\mathbb E}\left[\left(\partial_\theta v_\epsilon(t,\theta)
\right)^2\right]  \, = \, \langle Q^t
\rho_\epsilon'(\cdot-\theta), \rho_\epsilon'(\cdot-\theta)
\rangle.
\end{equation}
Therefore, by (\ref{cc}) and (\ref{ccc}):
\begin{eqnarray*}
& & \int_{\mathbb R} P_tG_{\epsilon,a}(z) \, \psi(a) \, da \, =
\int_0^1 d\theta \,
h_\theta  \int_{\mathbb R} {\cal N}(0,q_t(\theta)) (dy)
\, \psi(y+z(t,\theta))
\, \cdot
\\ \\ & & \quad \cdot 
\left[ \left(\partial_\theta z_\epsilon(t,\theta)\right)^2
- \, c^t_{\epsilon,\theta} + 2 \, y \,
\ell'_{\epsilon,\theta} \, \partial_\theta
z_\epsilon(t,\theta)  + \left(y^2-q_t(\theta)\right) 
\left(\ell'_{\epsilon,\theta} \right)^2 \right].
\end{eqnarray*}
Therefore we obtain:
\begin{eqnarray*}
& &
P_t G_\epsilon(z) \, = \, \int_0^1 h_\theta \,
\frac{e^{-(z(t,\theta))^2/2q_t(\theta)}} 
{\sqrt{2\pi q_t(\theta)}} \, \Big[ \left(\partial_\theta
z_\epsilon(t,\theta)\right)^2 - \, c^t_{\epsilon,\theta}
\\ \\ & & \qquad 
- 2 \, z(t,\theta) \,
\ell'_{\epsilon,\theta} \, \partial_\theta
z_\epsilon(t,\theta)  + \left(
(z(t,\theta))^2-q_t(\theta)\right) 
\left(\ell'_{\epsilon,\theta} \right)^2 \Big] \, d\theta,
\end{eqnarray*}
and
\[
\|P_t \, G_\epsilon\|^2 \, \leq \,
4\sum_{i=1}^3 I_i(t,\epsilon), \qquad I_i(t,\epsilon)
\, := \, \|V_{\epsilon,t}^i\|^2,
\]
where
\begin{eqnarray*}
V_{\epsilon,t}^1(z) & := & \int_0^1 h_\theta \,
\frac{e^{-(z(t,\theta))^2/2q_t(\theta)}} 
{\sqrt{2\pi q_t(\theta)}}
\left[ (\partial_\theta z_\epsilon(t,\theta))^2
-c^t_{\epsilon,\theta} \right] \, d\theta,
\\ \\
V_{\epsilon,t}^2(z) & := & - \int_0^1 h_\theta \,
\frac{e^{-(z(t,\theta))^2/2q_t(\theta)}} 
{\sqrt{2\pi q_t(\theta)}} \,
2 \, z(t,\theta) \,
\ell'_{\epsilon,\theta} \, \partial_\theta
z_\epsilon(t,\theta) \, d\theta,
\\ \\
V_{\epsilon,t}^3(z) & := & \int_0^1 h_\theta \, 
\frac{e^{-(z(t,\theta))^2/2q_t(\theta)}} 
{\sqrt{2\pi q_t(\theta)}} \, \left[
(z(t,\theta))^2 - q_t(\theta) \right]
\left(\ell'_{\epsilon,\theta} \right)^2 \, d\theta.
\end{eqnarray*}
For $F\in C^1_e(L)$, $k\in L$ and $K:=Qk\in H^1$
we have integrating by parts w.r.t. the Wiener measure: 
\[
{\mathbb E}\left[ \partial_K F(B) \right] \, = \,
{\mathbb E}\left[ F(B) \int_0^1 K_\theta' \,
dB_\theta\right]. 
\]
On the other hand, integrating by parts on $[0,1]$ we
obtain:
\[
\int_0^1 K_\theta' \, dB_\theta \, = \, K'_1 \, B_1 \, -
\, K'_0 \, B_0 \, - \int_0^1 K''_\theta \, B_\theta \,
d\theta \, = \, \int_0^1 k_\theta \, B_\theta \, d\theta,
\]
since $K'_1=B_0=0$. Therefore we obtain the following
formula: 
\begin{equation}\label{byi}
{\mathbb E}\left[F(B) \, \langle k,B\rangle\right] \, = \,
{\mathbb E}\left[ \partial_K F(B) \right].
\end{equation}
Iterating (\ref{byi}) several times we obtain for $F\in C^4_b(L)$,
$k^i\in L$ and $K^i:=Qk^i$:
\begin{eqnarray}\label{infe2}
& &
{\mathbb E}\left[F(B) \,
\langle k^1,B\rangle \, \langle k^2,B\rangle \right]
\\ \nonumber \\ \nonumber
& & = \ \langle K^1,k^2\rangle \,
{\mathbb E}\left[F(B)\right] \, + \, 
{\mathbb E}\left[\partial^2_{K^1,K^2} \, F(B) \right],
\end{eqnarray}
\begin{eqnarray}\label{infe4}
& &
{\mathbb E}\left[F(B) \, \langle k^1,B\rangle^2
\, \langle k^2,B\rangle^2 \right] 
\\ \nonumber \\ \nonumber
& & = \, \left(\langle K^1,k^1\rangle \,
\langle K^2,k^2\rangle \, + \, 2 \,
\langle K^1,k^2\rangle^2\right) {\mathbb E}\left[F(B) \right]
\\ \nonumber \\ \nonumber & &
\quad + \sum_{i\ne j}\langle K^i,k^i\rangle \,
{\mathbb E}\left[\partial^2_{K^j,K^j} \, F(B) \right]
\, + \, 4 \, \langle K^2,k^1\rangle \,
{\mathbb E}\left[\partial^2_{K^1,K^2} \, F(B) \right]
\\ \nonumber \\ \nonumber
& & \quad
+ \ {\mathbb E}\left[ \partial^4_{K^1,K^1,K^2,K^2}
\, F(B) \right].
\end{eqnarray}

\vspace{.3cm}\noindent
{\bf Estimate of $I_1$.} We set for the rest of the proof:
\begin{equation}\label{k^i}
k^1 := - \, e^{tA} \, \rho'_\epsilon(\cdot-\theta), \qquad
k^2 := - \, e^{tA} \, \rho'_\epsilon(\cdot-\theta'), \qquad
K^i \, := \, Qk^i,
\end{equation}
\[
F^{a,b}(z) \, := \,
\frac{e^{-(z(t,\theta)-a)^2/2q_t(\theta)}} 
{\sqrt{2\pi q_t(\theta)}} \, \cdot \,
\frac{e^{-(z(t,\theta')-b)^2/2q_t(\theta')}} 
{\sqrt{2\pi q_t(\theta')}}, \qquad
F \, := \, F^{0,0},
\]
for $z\in L$ and $a,b\in{\mathbb R}$. Then we have
\begin{eqnarray*}
& &
I_1(t,\epsilon) \, = \, 
\int \mu(dz) \left[ \int_0^1 h_\theta
\frac{e^{-(z(t,\theta))^2/2q_t(\theta)}} 
{\sqrt{2\pi q_t(\theta)}} \, \left[
(\partial_\theta z_\epsilon(t,\theta))^2 -
c^t_{\epsilon,\theta} \right] \, d\theta \right]^2 
\\ \\ & & = 
\int_{[0,1]^2} d\theta \, d\theta' \, h_\theta \,
h_{\theta'} \, \cdot
\\ \\ & & \cdot \ {\mathbb E}\left[ F(B) \, \left(
\langle k^1,B\rangle^2 \, \langle k^2,B\rangle^2 \, - \, 
\langle k^1,B\rangle^2 
\, c^t_{\epsilon,\theta'} \, - \, c^t_{\epsilon,\theta}
\, \langle k^2,B\rangle^2 \, + \, c^t_{\epsilon,\theta} \,
c^t_{\epsilon,\theta'}\right) \right].
\end{eqnarray*}
Moreover by (\ref{defce}) and (\ref{defcet}):
\[
\langle K^1,k^1\rangle \, = \, \langle Q^t
\rho'(\cdot-\theta), \rho'(\cdot-\theta) \rangle \, = \,
c^t_{\epsilon,\theta}, \qquad 
\langle K^2,k^2\rangle \, = \, c^t_{\epsilon,\theta'},
\]
\[
\langle K^1,k^2 \rangle \, = \,  \langle Q^t
\rho'(\cdot-\theta), \rho'(\cdot-\theta') \rangle \, =: \,
c^t_{\epsilon,\theta,\theta'}.
\]
Using (\ref{infe2}) and (\ref{infe4}), several terms
cancel and what remains is: 
\begin{eqnarray*}
& &
I_1(t,\epsilon) \, = \,
\int_{[0,1]^2} d\theta \, d\theta' \, h_\theta \,
h_{\theta'} \, \cdot 
\\ \\ & & \cdot \, {\mathbb E}\left[ F(B) \,
2 \, \langle K^1,k^2 \rangle^2 
+ \, 4 \, \langle K^1,k^2\rangle \,
\partial^2_{K^1,K^2} \, F(B) \, +
\ \partial^4_{K^1,K^1,K^2,K^2} \, F(B) \right].
\end{eqnarray*}
Notice that the function $\Gamma:{\mathbb
R}^2\mapsto{\mathbb R}_+$ 
\[
\Gamma(a,b) \, := \, {\mathbb E}\left[F^{a,b}(B)\right] 
\, = \, {\mathbb E}\left[ 
\frac{\exp\left(-\frac{(\langle
B,e^{tA}\delta_\theta\rangle-a)^2}{2q_t(\theta)}
\, - \, \frac{(\langle B,e^{tA}\delta_{\theta'}\rangle-b)^2}
{2q_t(\theta')} \right)} 
{2\pi \sqrt{q_t(\theta)\, q_t(\theta')}}
\right] 
\]
is the density of the convolution between ${\cal N}(0,q_t(\theta))
\otimes{\cal N}(0,q_t(\theta'))$ and the law of  
$(\langle B,e^{tA}\delta_\theta\rangle,
\langle B,e^{tA}\delta_{\theta'}\rangle)$. Therefore
$\Gamma$ is the density of the Gaussian measure on
${\mathbb R}^2$ with zero mean and covariance matrix:
\[
\begin{pmatrix} q_t(\theta) & 0 \\ \\ 0 & q_t(\theta') 
\end{pmatrix} \ + \ 
\begin{pmatrix} q^t(\theta) & q^t(\theta,\theta') 
\\ \\ q^t(\theta,\theta') &  q^t(\theta')
\end{pmatrix} \ = \ 
\begin{pmatrix} q_\infty(\theta) & q^t(\theta,\theta') 
\\ \\ q^t(\theta,\theta') & q_\infty(\theta') 
\end{pmatrix} =: \, \Lambda_{\theta,\theta'}. 
\]
Moreover
\begin{eqnarray*}
& & (q^t(\theta,\theta'))^2  \, = \, \left( {\mathbb
E}\left[ \langle B,e^{tA}\delta_\theta\rangle
\, \langle B,e^{tA}\delta_{\theta'}\rangle \right] \right)^2   
\\ \\ & & \leq \, {\mathbb
E}\left[ \langle B,e^{tA}\delta_\theta\rangle^2 \right]
{\mathbb E}\left[ \langle B,e^{tA}\delta_{\theta'}\rangle^2 
\right] \, = \, q^t(\theta) \, q^t(\theta')
\, \leq \, q^t(\theta) \, q_\infty(\theta').
\end{eqnarray*}
Using this inequality and recalling that $q_\infty-q^t=q_t$ 
we have:
\[
\det\Lambda_{\theta,\theta'} \, = \,
q_\infty(\theta) \, q_\infty(\theta') -
(q^t(\theta,\theta'))^2 \, \geq \,
q_t(\theta) \, q_\infty(\theta').
\]
Therefore by (\ref{estq}), for $\theta,\theta'\in[\delta,1-\delta]$:
\[
{\mathbb E}\left[F(B)\right] \, = \,
\Gamma(0,0) \, = \, \frac 1{2\pi(\det\Lambda_{\theta,\theta'})^{1/2}}
\, \leq \, \frac {\kappa_\delta^{-1/2}}{t^{1/4}}.
\]
Now, by (\ref{cos}):
\[
c^t_{\epsilon,\theta,\theta'} \, = \, \langle Q^t
\rho'_\epsilon(\cdot-\theta), \rho'_\epsilon(\cdot-\theta') \rangle 
\, = \, \sum_{i=1}^\infty
\frac{e^{-\lambda_i t}}{\lambda_i} \,
(\rho_\epsilon*e'_i)_\theta \,
(\rho_\epsilon*e'_i)_{\theta'}.  
\]
Setting $\eta_i:=\lambda_i^{1/2} e_i$ we have that
$(\eta_i)_{i\in{\mathbb N}}$ is a c.o.s. in $L$. We obtain
\begin{eqnarray*}
& &
\int_{[0,1]^2} d\theta \, d\theta' \, 
h_\theta \, h_{\theta'} \, {\mathbb
E}\left[F(B)\right] 
\left( c^t_{\epsilon,\theta,\theta'} \right)^2
\, \leq \, \frac \kappa{t^{1/4}}
\int_{[0,1]^2} d\theta \, d\theta' \, 
h_\theta \, h_{\theta'} \, 
\left( c^t_{\epsilon,\theta,\theta'} \right)^2 
\\ \\ & & = \, \frac \kappa{t^{1/4}}
\sum_{i,j=1}^\infty e^{-(\lambda_i+\lambda_j) t} \,
\left[ \int_0^1 (\rho_\epsilon*\eta_i)_\theta \,
(\rho_\epsilon*\eta_j)_{\theta}
\, h_\theta \,d\theta \right]^2.
\end{eqnarray*}
Now, since $\rho_\epsilon$ is a symmetric convolution
kernel:
\begin{eqnarray*}
& &
\int_0^1 (\rho_\epsilon*\eta_i)_\theta \,
(\rho_\epsilon*\eta_j)_{\theta}
\, h_\theta \,d\theta \, = \,
\langle\eta_j, \rho_\epsilon*[
h(\rho_\epsilon*\eta_i)] \rangle
\\ \\ & & \Longrightarrow \, 
\sum_{j=1}^\infty 
\left[ \int_0^1 (\rho_\epsilon*\eta_i)_\theta \,
(\rho_\epsilon*\eta_j)_{\theta}
\, h_\theta \,d\theta \right]^2 = \,
\|\rho_\epsilon*[h(\rho_\epsilon*\eta_i)]\|^2 \
\leq \, \|h\|^2,
\end{eqnarray*}
so that:
\[
\int_{[0,1]^2} d\theta \, d\theta' \, 
h_\theta \, h_{\theta'} \, {\mathbb
E}\left[F(B)\right] 
\left( c^t_{\epsilon,\theta,\theta'} \right)^2
\, \leq \, \frac{\kappa \, \|h\|^2}{t^{1/4}}
\sum_{i=1}^\infty e^{-\lambda_i 
t} \, \leq \, \frac{\kappa \, \|h\|^2}{t^{3/4}}.
\]
Now for all $\ell\in L$ we have:
\[
F(z+s\, \ell) \, = \,
F^{-s \, e^{tA}\ell_\theta,-s \, e^{tA}\ell_{\theta'}}(z)  \, = \,
F^{s \, e^{tA}\ell_\theta, \, s \, e^{tA}\ell_{\theta'}}(z) 
\]
so that, setting $H^i:=e^{tA}Qk^i$:
\begin{eqnarray}
& &  \label{2=0}
{\mathbb E}\left[ \partial^2_{K^1,K^2} \, F(B) \right]
\, = \,
\left. \frac{\partial^2}{\partial r \, \partial s}
{\mathbb E}\left[F(B+rK^1+sK^2) \right] \right|_{r=s=0}
\\ \nonumber \\ \nonumber & & 
= \, \left. \frac{\partial^2}{\partial r \, \partial s}
\, \Gamma(rH^1_\theta+sH^2_\theta, \,
rH^1_{\theta'}+sH^2_{\theta'})
\, \right|_{r=s=0} \, = \, - \, 
\frac{{\bf v}^T_1 \, \Lambda^{-1}_{\theta,\theta'} \,
{\bf v}_2} {2\pi\sqrt{\det\Lambda_{\theta,\theta'}}} \,
\end{eqnarray}
where ${\bf v}_i=(H^i_\theta,H^i_{\theta'})\in{\mathbb R}^2$. 
Since the entries of $\Lambda_{\theta,\theta'}$ are bounded 
uniformly 
in $\theta,\theta'\in[0,1]$ and for all $\theta,\theta'\in[0,1]$:
\begin{equation}\label{etaq}
| H^j_\theta | \, \leq \, \sum_{i=1}^\infty
\frac{e^{-\lambda_i t}} 
{\lambda_i} \, \|\rho_\epsilon*e_i'\|_\infty \,
\| e_i \|_\infty \, \leq \, 
\sum_{i=1}^\infty \frac{e^{-\lambda_i t}}{\lambda_i^{1/2}} 
\, \leq \, \kappa (1+|\ln t|),
\end{equation}
then we obtain:
\[
\left| {\mathbb E}\left[ \partial^2_{K^1,K^2} \, F(B)
\right] \right| \, \leq \,
\frac{\kappa(1+|\ln t|)^2}{(\det\Lambda_{\theta,\theta'})^{3/2}}
\]
and therefore:
\[
\int_{[0,1]^2} d\theta \, d\theta' \, h_\theta \,
h_{\theta'} \, \langle K^1,k^2\rangle \,
{\mathbb E}\left[ \partial^2_{K^1,K^2} \, F(B) \right]
 \, \leq \, \frac{\kappa(1+|\ln t|)^2}{t^{3/4}}.
\]
Analogously:
\begin{eqnarray*}
{\mathbb E}\left[ 
\partial^4_{K^1,K^1,K^1,K^2} \, F(B) \right] & = &
\left. \frac{\partial^4}{\partial^2 r \, \partial^2 s}
\, \Gamma(rH^1_\theta+sH^2_\theta, \,
rH^1_{\theta'}+sH^2_{\theta'})
\, \right|_{r=s=0} 
\\ \\ & = & \frac 1{(\det\Lambda_{\theta,\theta'})^{3/2}} 
\, R_{\theta,\theta'}
(H^1_\theta,H^2_\theta,H^1_{\theta'},H^2_{\theta'}),
\end{eqnarray*}
where $R_{\theta,\theta'}$ is a multi-linear form on ${\mathbb 
R}^4$ with uniformly bounded coefficients w.r.t.
$\theta,\theta'\in[0,1]$. Therefore:
\[
\int_{[0,1]^2} d\theta \, d\theta' \, h_\theta \,
h_{\theta'} \, {\mathbb E}\left[
\partial^4_{K^1,K^1,K^1,K^2} \, F(B)\right]
\, \leq \, \frac{\kappa(1+|\ln t|)^4}{t^{3/4}}.
\]

\vspace{.3cm}\noindent
{\bf Estimate of $I_2$.} Continuing with the notations
introduced in the previous step, we notice now that:
\[
\frac{e^{-(z(t,\theta))^2/2q_t(\theta)}} 
{\sqrt{2\pi q_t(\theta)}} \cdot
\frac{z(t,\theta)}{q_t(\theta)} \cdot
\frac{e^{-(z(t,\theta))^2/2q_t(\theta)}} 
{\sqrt{2\pi q_t(\theta)}} \cdot
\frac{z(t,\theta)}{q_t(\theta)} \, = \,
\left. \frac{\partial^2}{\partial a \, \partial b}
\, F^{a,b}(z) \right|_{a=b=0}.
\]
Then, setting $\nu_{\epsilon,\theta} :=
(\rho_\epsilon* q_t(\cdot,\theta))_\theta'$, we have:
\begin{eqnarray*}
& &
I_2(t,\epsilon) \, = \,
\int \mu(dz) \left[ \int_0^1 h_\theta \,
\frac{e^{-(z(t,\theta))^2/2q_t(\theta)}} 
{\sqrt{2\pi q_t(\theta)}} \,
2 \, \frac{z(t,\theta)}{q_t(\theta)} \,
\nu_{\epsilon,\theta} \, \partial_\theta
z_\epsilon(t,\theta) \, d\theta \right]^2 
\\ \\ & & = \, 4
\int_{[0,1]^2} d\theta \, d\theta' \, h_\theta \,
h_{\theta'} \, \nu_{\epsilon,\theta}
\, \nu_{\epsilon,\theta'} \ {\mathbb E}\left[
\langle B, k^1 \rangle \, \langle B, k^2 \rangle \,
\left. \frac{\partial^2}{\partial a \, \partial b}
\, F^{a,b}(B) \right|_{a=b=0} \right].
\end{eqnarray*}
By (\ref{infe2}) we have:
\begin{eqnarray*}
& &
{\mathbb E}\left[
\langle B, k^1 \rangle \, \langle B, k^2 \rangle \, 
F^{a,b}(B) \right]
\\ \\ & &   = \, \langle K^1,k^2 \rangle 
\, {\mathbb E}\left[F^{a,b}(B)\right]
+ {\mathbb E}\left[\partial^2_{K^1,K^2} \, F^{a,b}(B)\right] 
\\ \\ & & = \, c^t_{\epsilon,\theta,\theta'} \, \Gamma(a,b) 
\, + \, \left. \frac{\partial^2}{\partial r \, \partial s}
\, \Gamma(a-rH^1_\theta-sH^2_\theta, \,
b-rH^1_{\theta'}-sH^2_{\theta'}) \, \right|_{r=s=0}.
\end{eqnarray*}
Now, recalling that $\Gamma$ is the density of 
${\cal N}(0,\Lambda_{\theta,\theta'})$, we can compute:
\[
\left. \frac{\partial^2}{\partial a \, \partial b}
\, \Gamma(a,b) \right|_{a=b=0} \, = \
\frac{q^t(\theta,\theta')}
{2\pi(\det\Lambda_{\theta,\theta'})^{3/2}},
\]
\begin{eqnarray*}
& &
\left. \frac{\partial^4}{\partial a \, \partial b \,
\partial r \, \partial s} \Gamma(a-rH^1_\theta-sH^2_\theta, \,
b-rH^1_{\theta'}-sH^2_{\theta'})
\,  \right|_{a=b=r=s=0} 
\\ \\ & & = \, \frac 1{(\det\Lambda_{\theta,\theta'})^{3/2}} 
\, \widehat R_{\theta,\theta'}
(1,1,H^1_\theta,H^2_\theta,H^1_{\theta'},H^2_{\theta'})
\end{eqnarray*}
where $\widehat R_{\theta,\theta'}$ is a multi-linear form on ${\mathbb 
R}^6$ with uniformly bounded coefficients w.r.t.
$\theta,\theta'\in[0,1]$. By (\ref{cos}):
\[
q_t(\theta,\theta') \, = \,  \theta\wedge\theta' - 
\sum_{i=1}^\infty \frac{e^{-\lambda_i t}} 
{\lambda_i} \, e_i(\theta) \ e_i(\theta'). 
\]
Since $(\rho_\epsilon* q_\infty(\cdot,\theta))_\theta'
=(\rho_\epsilon*1_{[0,\theta]})_\theta=1/2$, then:
\[
\nu_{\epsilon,\theta} \, = \,
(\rho_\epsilon* q_t(\cdot,\theta))_\theta'
= \, \frac 12 \, -
\sum_{i=1}^\infty \frac{e^{-\lambda_i t}} 
{\lambda_i} \, (\rho_\epsilon*e'_i)_\theta \ e_i(\theta)
\]
and therefore
\[
|\nu_{\epsilon,\theta}| \, \leq \, \kappa(1+|\ln t|). 
\]
Therefore we have proven that:
\begin{eqnarray*}
I_2(t,\epsilon) & \leq & \kappa
\int_{[0,1]^2} d\theta \, d\theta' \, h_\theta \,
h_{\theta'} \, \nu_{\epsilon,\theta} \, \nu_{\epsilon,\theta'} 
\, \frac{1+|\widehat R_{\theta,\theta'}
(1,1,H^1_\theta,H^2_\theta,H^1_{\theta'},H^2_{\theta'})|}
{(\det\Lambda_{\theta,\theta'})^{3/2}}
\\ \\ & \leq & \frac{\kappa(1+|\ln t|)^6}{t^{3/4}}.
\end{eqnarray*}

\vspace{.3cm}\noindent
{\bf Estimate of $I_3$.} Arguing like for $I_2$ we obtain:
\begin{eqnarray*}
I_3(t,\epsilon) & = &
\int \mu(dz) \left[ \int_0^1 h_\theta \,
\frac{e^{-(z(t,\theta))^2/2q_t(\theta)}} 
{\sqrt{2\pi q_t(\theta)}} \left[
\frac{(z(t,\theta))^2}{(q_t(\theta))^2}
\, - \, \frac 1{q_t(\theta)} \right]
\nu_{\epsilon,\theta}^2 \, d\theta \right]^2 
\\ \\ & = & \int_{[0,1]^2}
d\theta \, d\theta' \, h_\theta \, h_{\theta'} 
\, \nu_{\epsilon,\theta}^2 \, \nu_{\epsilon,\theta'}^2  \, 
\left. \frac{\partial^4}{\partial^2 a \, \partial^2 b} \,
\Gamma(a,b)\right|_{a=b=0}
\, \leq \, \frac{\kappa(1+|\ln t|)^6}{t^{3/4}},
\end{eqnarray*}
and the proof of Lemma \ref{convu} is complete. \quad
$\square$

\medskip\noindent
Using the proofs of Proposition \ref{pr3} and Lemma
\ref{convu}, we prove also the following:
\begin{corollary}\label{dir}
For all $h\in C_c(0,1)$, ${\cal R}G_\epsilon$ converges
weakly in ${\rm Dom}({\cal E})$ to ${\cal R}G_0\in {\rm Dom}({\cal E})$,
where for $\mu$-a.e $z\in C$:
\begin{eqnarray*}
& &
{\cal R}G_0(z) \, := \, \int_0^\infty \int_0^1 h_\theta \, 
\frac{e^{-(z(t,\theta))^2/2q_t(\theta)}} 
{\sqrt{2\pi q_t(\theta)}} \, \Bigg[ \left(
\frac{\partial z(t,\theta)}{\partial\theta}\right)^2
- \, c^t_{0,\theta}
\\ \\ & & \qquad 
- 2 \, \nu_{0,\theta} \, \frac{z(t,\theta)}{q_t(\theta)} \, 
\frac{\partial z(t,\theta)}{\partial\theta} \,  + \nu_{0,\theta}^2 
\left( \left[\frac{z(t,\theta)}{q_t(\theta)} \right]^2-
\frac 1{q_t(\theta)}\right) 
\Bigg] \, d\theta \, dt, 
\end{eqnarray*}
for $\theta\in(0,1)$, $t\in(0,\infty)$, $z(t,\theta)$ is defined
by (\ref{defzv}) and:
\[
c^t_{0,\theta} \, := \, \sum_{i=1}^\infty
\frac{e^{-\lambda_i t}}{\lambda_i} \,
\left|e'_i(\theta)\right|^2, \qquad
\nu_{0,\theta} \, := \, \frac 12 \, - \sum_{i=1}^\infty 
\frac{e^{-\lambda_i t}}{\lambda_i} \,
e'_i(\theta) \, e_i(\theta).
\]
Moreover for all $F\in {\rm Lip}_e(L)$ and $a\in{\mathbb R}$:
\begin{eqnarray}\label{expl}
& &
{\mathbb E}\left[ F(B) \int_0^1 h_\theta 
\, : \! \B_\theta^2 \! : \, dL^a_\theta \right] 
\\ \nonumber \\ \nonumber & & = \,
{\mathbb E}[ F(B) ] \, \int_0^1 h_\theta \, 
\frac{a^2-\theta}{4\theta^2} \,
\frac{e^{-a^2/2\theta}}{\sqrt{2\pi\theta}} \, d\theta \, + \,
{\cal E} \!  \left(e^{a\langle l'', \cdot \rangle} \ F, \
{\cal R}G_0\right) \ e^{-a^2\|\ell'\|^2/2},
\end{eqnarray}
where $l\in C^2([0,1])$, $l(0)=0$ and $l(x)=1$ for all $x$
such that $h(x)\ne 0$.  
\end{corollary}
\noindent 
Formula (\ref{expl}) allows to compute directly the value of the generalized
functional constructed in Theorem \ref{main1} without using the limit in the
l.h.s. of (\ref{lim}).

\section{The case of quadratic $\varphi$ and constant $h$}

We want to consider the divergence of a vector field
of particular interest, namely the identity ${\cal
K}(\omega) = \omega$. This case corresponds to
$\varphi(r)=\frac 12 r^2$ and $h\equiv 1$, and therefore it
does not fit in the assumptions of Theorem \ref{main1.5},
since $h$ has not compact support in $(0,1)$. Still, since
$\varphi''\equiv 1$, this
case is simpler than the general one and can be treated
without the main estimate of Lemma \ref{convu}.

Let us go back to the result of Lemma \ref{u}: formula
(\ref{u2}) becomes 
\[
k_\theta \, {\mathbb E}[ B_\theta+K_\theta ] \, = \,
- \, \frac 12 \, \frac{d^2}{d\theta^2} \, {\mathbb E}\left[
(B_\theta+K_\theta)^2 \right]
+ \, {\mathbb E}\left[ \lambda(\theta, K'_{\theta}, B_\theta)
\right],
\]
\[
{\rm i.e.} \qquad
k_\theta \, K_\theta \, = \,
- \, \frac 12 \, \frac{d^2}{d\theta^2} \, 
(\theta+K_\theta^2) \, + \, (K'_{\theta})^2.
\]
Integrating over $[0,1]$ in $d\theta$ we obtain:
\[
\int_0^1 k_\theta \, K_\theta \, d\theta \, = \,
- \frac 12 \left[ (K_\theta^2)' \right]_0^1 \, + \,
\int_0^1 (K'_{\theta})^2 \, d\theta \, = \,
\int_0^1 (K'_{\theta})^2 \, d\theta,
\]
since $K_0=K_1'=0$. By (\ref{cme}) this yields for $\Psi_k :=
e^{\langle k,\cdot\rangle}$:
\begin{eqnarray*}
{\mathbb E}\left[ \partial_B \Psi_k(B) \right] & = &
e^{\frac 12\langle Qk,k\rangle}
\int_0^1 (K'_{\theta})^2 \, d\theta \, = \,
{\mathbb E}\left[ \Psi_k(B) \int_0^1 :\! \B_\theta^2 \! :
\, d\theta \right]
\\ \\ & := &
\lim_{\epsilon\to 0} 
{\mathbb E}\left[ \Psi_k(B) \int_0^1 :\! 
\B_{\epsilon,\theta}^2 \! : \, d\theta \right].
\end{eqnarray*}
In this case $:\! \B_\theta^2 \! :$ 
appears without integration w.r.t. the local time
process and is therefore defined in the classical way, see
\cite{hkps}. Arguing like in sections 5 and 6, we set now:
\[
{\cal G}_\epsilon(B) \, := \, \int_0^1 :\! 
\B_{\epsilon,\theta}^2 \! : \, d\theta,
\]
and we compute for all $z\in C$:
\[
P_t {\cal G}_\epsilon (z) \, = \, \int_0^1 \left[
\left(\partial_\theta z_\epsilon(t,\theta) \right)^2
- c_{\epsilon,\theta}\right] \, d\theta.
\]
Arguing like in the proof of Lemma \ref{convu}, see in
particular the estimate of $I_1$, we compute:
\begin{eqnarray*}
& & \|P_t {\cal G}_\epsilon\|^2_{L^2(\mu)} \, = \, 
\int \mu(dz) \left[ \int_0^1 \left[
(\partial_\theta z_\epsilon(t,\theta))^2 -
c^t_{\epsilon,\theta} \right] \, d\theta \right]^2 
= \int_{[0,1]^2} d\theta \, d\theta' \, \cdot 
\\ \\ & & \cdot \, {\mathbb E}\left[ 
\langle k^1,B\rangle^2 \, \langle k^2,B\rangle^2 \, - \, 
\langle k^1,B\rangle^2 
\, c^t_{\epsilon,\theta'} \, - \, c^t_{\epsilon,\theta}
\, \langle k^2,B\rangle^2 \, + \, c^t_{\epsilon,\theta} \,
c^t_{\epsilon,\theta'} \right]
\\ \\ & & = \int_{[0,1]^2} d\theta \, d\theta' \,
2 \, \langle Qk^1,k^2\rangle^2 \, = 2 \int_{[0,1]^2} d\theta \, 
d\theta' \left[\sum_{i=1}^\infty
\frac{e^{-\lambda_i t}}{\lambda_i} \,
(\rho_\epsilon*e'_i)_\theta \,
(\rho_\epsilon*e'_i)_{\theta'}\right]^2
\\ \\ & & = 2 \sum_{i,j=1}^\infty e^{-(\lambda_i+\lambda_j) t}
\left[ \int_0^1 (\rho_\epsilon*\eta_i)_\theta \,
(\rho_\epsilon*\eta_j)_{\theta} \,d\theta \right]^2 
\leq \, 2\sum_{i=1}^\infty e^{-\lambda_i t} \,
\leq \, \frac\kappa{t^{1/2}}.
\end{eqnarray*}
Therefore ${\cal R}{\cal G}_\epsilon$ converges weakly in 
${\rm Dom}({\cal E})$ to ${\cal R}{\cal G}_0 \in {\rm Dom}({\cal E})$ and
\[
{\mathbb E}[ \langle \nabla F(B), B\rangle ] \, = \,
\lim_{\epsilon\to 0} \,
{\mathbb E}\left[ F(B) \int_0^1 :\! 
\B_{\epsilon,\theta}^2 \! : \, d\theta \right] \, = \,
{\cal E}(F,{\cal R}{\cal G}_0), 
\]
for all $F\in {\rm Dom}({\cal E})$, where:
\[
{\cal R}{\cal G}_0 (z) \, = \, \int_0^\infty
\left[ \left(
\frac{\partial z(t,\theta)}{\partial\theta}\right)^2
- \, c^t_{0,\theta} \right] dt, \qquad \mu-{\rm a.e.} \ z\in C,
\]
see Corollary \ref{dir} above.

\end{document}